\documentclass[11pt]{amsart}
\usepackage[utf8]{inputenc}
\usepackage[T1]{fontenc}
\usepackage{lmodern}
\usepackage[margin=1.05in]{geometry}
\usepackage{amsmath,amssymb,amsthm}
\usepackage{enumerate}
\usepackage[colorlinks=true,linkcolor=blue,citecolor=blue,urlcolor=blue]{hyperref}

\allowdisplaybreaks

\newtheorem{theorem}{Theorem}[section]
\newtheorem{lemma}[theorem]{Lemma}
\newtheorem{proposition}[theorem]{Proposition}
\newtheorem{corollary}[theorem]{Corollary}
\theoremstyle{definition}
\newtheorem{definition}[theorem]{Definition}
\theoremstyle{remark}
\newtheorem{remark}[theorem]{Remark}
\numberwithin{equation}{section}

\newcommand{\R}{\mathbb{R}}
\newcommand{\Rn}{\mathbb{R}^n}
\newcommand{\Z}{\mathbb{Z}}
\newcommand{\abs}[1]{\left|#1\right|}

\newcommand{\B}[3]{\dot B^{#1}_{#2,#3}}
\newcommand{\Dj}[1]{\dot\Delta_{#1}}
\newcommand{\Sj}[1]{\dot S_{#1}}
\newcommand{\tDj}[1]{\widetilde{\dot\Delta}_{#1}}
\newcommand{\F}{\mathcal{F}}
\newcommand{\Sch}{\mathcal{S}}
\newcommand{\dK}{\dot K}
\newcommand{\Ic}{\mathcal{I}}
\newcommand{\E}{E_{s,\zeta}}
\newcommand{\supp}{\operatorname{supp}}

\begin{document}

\title[Hardy--Hénon wave equation ]{Global well-posedness and asymptotic behavior for the Hardy--Hénon wave equation}

\author{Edison Cuba}
\address{Departamento Académico de Matemáticas y Estadística, Universidad Nacional de San
Antonio Abad del Cusco,  Av. de la Cultura 733, Cusco 08003, Peru}
\email{edison.cuba@unsaac.edu.pe}

\subjclass[2020]{35L71, 35L05, 35C06, 35B40, 46E35}
\keywords{Hardy--Hénon wave equation, Global
well-posedness, composition estimate,  Besov-type space}

\begin{abstract}
We study the Cauchy problem for the Hardy--Hénon wave equation, a semilinear wave equation in which a real power nonlinearity is multiplied by the singular weight $|x|^{-\gamma}$, and carry out the analysis in higher dimensions, $n\geq 2$. In a time-weighted framework built on homogeneous Besov space, we prove global well-posedness for small data. The class of initial data contains critical Besov spaces and homogeneous functions, which yields self-similar solutions. We also prove scattering, decay, asymptotic stability and asymptotic self-similarity. The proofs combine dispersive estimates for the wave group, the product estimates to remove the weight, and new composition and difference estimates for the fractional power. 
\end{abstract}

\maketitle

\setcounter{tocdepth}{2}

\section{Introduction}

Semilinear evolution equations in which a power of $|x|$ multiplies the nonlinear term have been studied intensively in recent years. The name Hardy--Hénon refers to the work of Hénon \cite{Hen73}, where the elliptic counterpart of these equations arose as a model in astrophysics, and to the Hardy inequality, which governs the case of singular coefficients. In this paper we study the Cauchy problem for the Hardy--Hénon wave equation
\begin{equation}
\left\{
\begin{array}{l}
\partial _{t}^{2}u\,-\,\Delta u\,=\,a\,\abs x^{-\gamma}\abs u^{p-1}u, \\
u(0)=u_{0},\quad \partial _{t}u(0)=u_{1},%
\end{array}%
\right. \label{eq:HHW}
\end{equation}%
in $\Rn$ with $n\ge2$, where $a\in\R$, $0<\gamma<\min\{2,n\}$, and the power $p>1$ is a real number. The unknown $u$ is real valued, and the sign of $a$ plays no role in our arguments.

Since the coefficient $|x|^{-\gamma}$ is singular at the origin, classical solutions cannot be expected in general, and we work with mild solutions, that is, solutions of the Duhamel formulation of \eqref{eq:HHW}. The equation is invariant under a one-parameter group of dilations, and the solutions that are invariant under these dilations are called self-similar. Their initial data are homogeneous functions of degrees $-\frac{2-\gamma}{p-1}$ and $-\frac{2-\gamma}{p-1}-1$, which belong to no Sobolev space $\dot H^s$. Self-similar solutions are therefore out of reach of energy methods. They require function spaces that contain homogeneous functions, such as the homogeneous Besov spaces with third index $\infty$ used in this paper.

For $\gamma=0$, \eqref{eq:HHW} is the classical semilinear wave equation, and small data global solutions and self-similar solutions have been constructed in several settings. Pecher \cite{Pec00} used a time-weighted framework based on Besov spaces with third index $2$, Ribaud and Youssfi \cite{RY02} worked in homogeneous Sobolev spaces, and Planchon \cite{Pla00} treated integer powers in $L^2$-based Besov spaces with third index $\infty$. In \cite{FPL20}, the authors introduced, for integer powers, a time-weighted framework built on a single homogeneous Besov space with third index $\infty$. Cuba and Ferreira \cite{CF25} extended this framework to real powers and to weak-Besov spaces, by means of composition estimates for the fractional power. For wave equations with singular potentials, a related time-weighted approach in Lorentz spaces was developed in \cite{dAF17}. The aim of the present paper is to extend this framework to the singular nonlinearity $|x|^{-\gamma}|u|^{p-1}u$.

Our first result, Theorem \ref{thm:global}, provides global well-posedness for small data in a time-weighted class of initial data, under explicit conditions on $n$, $\gamma$, $p$ and the Lebesgue exponent $r$ of the solution space. Solutions with homogeneous data are self-similar. The class of data contains critical Besov spaces and smooth homogeneous functions (Corollary \ref{cor:data}), and the conditions can be met in every dimension $n\ge2$ (Remark \ref{rem:hypotheses}). We then study the behavior of the solutions for large times, following the approach developed by  \cite{BFP19} for the two-power Schrödinger equation. We prove scattering (Theorem \ref{thm:scattering}), faster decay for data whose free evolution decays faster (Theorem \ref{thm:decay}), asymptotic stability (Theorem \ref{thm:stability}), and asymptotic self-similarity for small perturbations of homogeneous data (Corollary \ref{cor:asymptotic-ss}). To the best of our knowledge, these are the first results on self-similar and asymptotically self-similar solutions of \eqref{eq:HHW} with $\gamma>0$.

The proofs combine three ingredients. The first is a dispersive estimate for the wave group, which maps Besov spaces built on $L^{r'}$ into Besov spaces built on $L^r$, with a decay rate in time and a loss of derivatives (Lemma \ref{lem:group}). The second ingredient removes the singular weight. The function $|x|^{-\gamma}$ belongs to homogeneous Besov spaces with third index $\infty$ (Proposition \ref{prop:weight}), and the product law of Abidi \cite{Abi07} applies to factors of this type. This gives a multiplier estimate for $|x|^{-\gamma}$ (Lemma \ref{lem:potential}) and a dispersive estimate involving the weight (Lemma \ref{lem:group-potential}), which is the wave analog of a heat estimate of \cite{Chi19}. The third ingredient, which is the main technical novelty, concerns the fractional power. Once the weight is removed, $|u|^{p-1}u$ has to be estimated in a Besov space of positive regularity. Composition estimates, which go back to an idea of  \cite{MYZ07}, are designed for the negative regularity that appears for the heat equation. We refine  them with a lemma on the high frequencies of the coefficients of a telescopic expansion of the nonlinearity (Lemma \ref{lem:gain}). It yields a gain of derivatives that is limited only by the smoothness of the nonlinearity. As a result, the regularity of $|u|^{p-1}u$ in the composition and difference estimates may be any number below $p-1$ and $p-2$, respectively (Lemmas \ref{lem:composition} and \ref{lem:difference}). Such gains are indispensable in high dimensions, because the regularity at which the nonlinearity is measured grows like $\frac n2$.

\subsection{Notation}

We assume $n\ge2$ throughout, unless stated otherwise. For $1\le r\le\infty$, the conjugate exponent $r'$ is defined by $\frac1r+\frac1{r'}=1$. For $0<\rho\le\infty$, $L^\rho(\Rn)$ is the usual Lebesgue space, which is only quasi-normed when $\rho<1$, and we write $\|\cdot\|_{L^\rho}$ for its (quasi-)norm. The space $\ell^\rho(\Z)$ consists of the sequences with summable $\rho$-th powers. The Schwartz class is $\Sch(\Rn)$ and its dual is $\Sch'(\Rn)$. We write $\F u=\widehat u$ and $\F^{-1}u$ for the Fourier transform and its inverse. We also set
\[
h_\gamma(x):=\abs x^{-\gamma},\qquad F(u):=\abs u^{p-1}u .
\]
The letter $C$ denotes a positive constant whose value may change from line to line, and $A\lesssim B$ means $A\le CB$.

\begin{definition}[{\cite[p.~22]{BCD11}}]\label{def:Sh}
We denote by $\Sch'_h(\Rn)$ the space of tempered distributions $u$ such that $\lim_{\lambda\to\infty}\|\F^{-1}[\theta(\lambda\,\cdot)\widehat u]\|_{L^\infty}=0$ for every $\theta\in C^\infty_c(\Rn)$.
\end{definition}

We fix a homogeneous Littlewood--Paley decomposition $\{\phi_j\}_{j\in\Z}$. Let $\phi\in C^\infty_c(\Rn)$ be non-negative and radial, with
\[
\supp\phi\subset\{\xi\in\Rn:\ 2^{-1}\le|\xi|\le2\}\quad\text{and}\quad\sum_{j\in\Z}\phi(2^{-j}\xi)=1\ \text{ for all }\xi\ne0.
\]
For $j\in\Z$ we put $\phi_j(\xi):=\phi(2^{-j}\xi)$ and $\Phi_j(\xi):=1-\sum_{k\ge j+1}\phi_k(\xi)$. Then
\[
\supp\Phi_j\subset\{|\xi|\le2^{j+1}\}\quad\text{and}\quad\supp\phi_j\subset\{2^{j-1}\le|\xi|\le2^{j+1}\}.
\]
The localization and cut-off operators are
\begin{equation}\label{eq:LP}
\Dj ju:=\F^{-1}[\phi_j\widehat u],\qquad \tDj ju:=\sum_{|k-j|\le2}\Dj ku,\qquad \Sj ju:=\F^{-1}[\Phi_j\widehat u],\qquad j\in\Z.
\end{equation}
The operators $\Sj j$ act on all of $\Sch'(\Rn)$, and $\Sj ju=\sum_{j'\le j}\Dj{j'}u$ whenever $u\in\Sch'_h(\Rn)$. Moreover, $\Dj j\tDj j=\Dj j$.

\begin{definition}[Homogeneous Besov space]\label{def:Besov}
For $s\in\R$ and $1\le r,\sigma\le\infty$, the homogeneous Besov space $\B s r\sigma(\Rn)$ is
\[
\B s r\sigma(\Rn):=\bigl\{u\in\Sch'_h(\Rn):\ \|u\|_{\B s r\sigma}<\infty\bigr\},\qquad
\|u\|_{\B s r\sigma}:=\Bigl\|\bigl(2^{js}\|\Dj ju\|_{L^r}\bigr)_{j\in\Z}\Bigr\|_{\ell^\sigma}.
\]
\end{definition}

This space is complete if $s<n/r$, or if $s=n/r$ and $\sigma=1$ (see \cite[Chapter 2]{BCD11}).

\subsection{Functional framework}

We now make the setting precise. Equation \eqref{eq:HHW} is invariant under the scaling
\begin{equation}\label{eq:scaling}
u_\lambda(t,x):=\lambda^{\frac{2-\gamma}{p-1}}u(\lambda t,\lambda x),\qquad \lambda>0,
\end{equation}
which transforms the data into $(\lambda^{\frac{2-\gamma}{p-1}}u_0(\lambda\,\cdot),\lambda^{\frac{2-\gamma}{p-1}+1}u_1(\lambda\,\cdot))$. A solution is self-similar when $u_\lambda=u$ for every $\lambda>0$. Self-similar solutions must have homogeneous data, namely
\begin{equation}\label{eq:homogeneous}
u_0(x)=\Omega_0\Bigl(\frac x{|x|}\Bigr)|x|^{-\frac{2-\gamma}{p-1}},\qquad u_1(x)=\Omega_1\Bigl(\frac x{|x|}\Bigr)|x|^{-\frac{2-\gamma}{p-1}-1},
\end{equation}
with functions $\Omega_0,\Omega_1$ defined on the unit sphere $\mathbb S^{n-1}$. Let $K(t):=\frac{\sin(t|D|)}{|D|}$ and $\dK(t):=\cos(t|D|)$, where $|D|:=(-\Delta)^{\frac12}$. For $t>0$ we look for a function (or a distribution) $u$ such that
\begin{equation}\label{eq:integral}
\begin{aligned}
&u(t)=u_L(t)+B(u)(t),\qquad u_L(t):=\dK(t)u_0+K(t)u_1,\\
&B(u)(t):=a\int_0^tK(t-\tau)\bigl(h_\gamma F(u(\tau))\bigr)\,d\tau .
\end{aligned}
\end{equation}
Solutions of \eqref{eq:integral} are called mild solutions of \eqref{eq:HHW}. The equation is time reversible, so it is enough to consider positive times.

For $2\le r\le\infty$ we introduce the numbers $\beta_r$ and $\eta_r$ through
\begin{equation}\label{eq:beta-eta}
\frac{2\beta_r}{n+1}=\frac{\eta_r}{n-1}=\frac12-\frac1r .
\end{equation}
The number $\eta_r$ is the rate of decay of the wave group from $L^{r'}$ to $L^r$, and $2\beta_r$ is the number of derivatives lost in this estimate (Lemma \ref{lem:group}). Note that $\eta_r<2\beta_r$ when $r>2$. For $p>1$, $0<\gamma<\min\{2,n\}$ and $2<r<\infty$ we set
\begin{equation}\label{eq:indices}
s_r:=\frac nr-\frac{1-\gamma+\eta_r}{p-1},\qquad
\zeta_r:=\frac{1-\eta_r}{p-1},\qquad
\sigma_r:=s_r-1+2\beta_r+\gamma .
\end{equation}
The index $s_r$ is chosen so that $F(u)$, measured in $\B{\sigma_r}{r'}{\infty}$, scales like $\|u\|_{\B{s_r}r\infty}^p$. In other words,
\begin{equation}\label{eq:homogeneity}
\sigma_r-\frac n{r'}=p\Bigl(s_r-\frac nr\Bigr).
\end{equation}
The scaling \eqref{eq:scaling} then determines $\zeta_r$, and
\begin{equation}\label{eq:zeta}
\zeta_r=\frac{2-\gamma}{p-1}+s_r-\frac nr,\qquad (p-1)\zeta_r+\eta_r=1 .
\end{equation}
Both \eqref{eq:homogeneity} and \eqref{eq:zeta} follow from \eqref{eq:beta-eta} and \eqref{eq:indices} by elementary algebra. In particular, the quantity $\sup_{t>0}t^{\zeta_r}\|u(t)\|_{\B{s_r}r\infty}$ is invariant under \eqref{eq:scaling}, up to a multiplicative constant. For $\gamma=0$, the numbers $s_r$ and $\zeta_r$ are the indices of the time-weighted framework of the classical wave equation, see Remark \ref{rem:hypotheses} (a).

\begin{definition}\label{def:E}
Let $s\in\R$, $\zeta>0$ and $2\le r\le\infty$. Let $\E$ be the set of maps $u:(0,\infty)\to\B sr\infty(\Rn)$ such that, for every $j\in\Z$, the function $(t,x)\mapsto(\Dj ju(t))(x)$ is measurable on $(0,\infty)\times\Rn$, and
\[
\|u\|_{\E}:=\sup_{t>0}t^\zeta\|u(t)\|_{\B sr\infty}<\infty .
\]
With $u_L$ as in \eqref{eq:integral}, define the class of initial data
\[
\Ic=\Ic(s,\zeta):=\bigl\{(u_0,u_1)\in\Sch'_h(\Rn)\times\Sch'_h(\Rn):\ u_L\in\E\bigr\},\qquad
\|(u_0,u_1)\|_{\Ic}:=\sup_{t>0}t^\zeta\|u_L(t)\|_{\B sr\infty}.
\]
\end{definition}

\subsection{Comparison with previous results}\label{subsec:comparison}

We first compare Theorem \ref{thm:global} with recent works on wave equations with a singular coefficient in the nonlinear term. In \cite{Pha24} the equation $\partial_t^2u-\Delta u+c_1|x|^{-2}u=c_2|x|^{-b}F(u)$ is considered for $0<b<2$ and $n\ge5$, with a nonlinearity $F$ of power type. For small radial data and small $c_1$, global solutions are obtained in weak-$L^\rho$ spaces without time weight, together with a scattering property and polynomial stability. Jiang, Shen and Li \cite{JSL26} proved local and global well-posedness of $\partial_t^2u-\Delta u+|x|^{-b}|u|^{\alpha}u=0$ in $\dot H^1(\R^3)\times L^2(\R^3)$ and in $\dot H^{s+1}(\R^3)\times\dot H^s(\R^3)$, in the energy subcritical regime, by means of Strichartz estimates. Lucente \cite{Luc00} studied $\partial_t^2u-\Delta u+V(x)|u|^{p-1}u=0$ with a positive and regular coefficient $V$, using weighted nonlinear estimates in Besov spaces. These results and ours do not contain each other. For instance, the linear part in \cite{Pha24} contains a Hardy potential, the data in \cite{JSL26} have finite energy, and the coefficient in \cite{Luc00} is regular. Theorem \ref{thm:global} and Corollary \ref{cor:data} allow homogeneous data, which lie in no space $\dot H^s$, impose no symmetry, and hold in every dimension $n\ge2$. To the best of our knowledge, self-similar and asymptotically self-similar solutions of \eqref{eq:HHW} with $\gamma>0$ have not been constructed before.

The results of Section \ref{sec:asymptotic} are the wave counterparts of known results for the Hardy--Hénon parabolic equation. Ben Slimene, Tayachi and Weissler \cite{BSTW17} constructed global self-similar and asymptotically self-similar solutions in Lebesgue spaces. In \cite{CIT22} was obtained optimal well-posedness and forward self-similar solutions in critical weighted Lebesgue spaces. In \cite{CITT25}, the authors constructed global solutions with data in weighted Lorentz spaces that behave asymptotically like self-similar solutions, and they proved that the asymptotic profiles are stable. For $\gamma=0$, asymptotically self-similar solutions of the wave equation were obtained by Pecher \cite{Pec00} in a framework based on $\dot B^s_{r,2}$.

We now turn to the nonlinear estimates. Lemmas \ref{lem:composition} and \ref{lem:difference} contain \cite[Lemmas 3.2 and 3.7]{Chi19} as the case $\nu=0$ (Remarks \ref{rem:Chikami-recovered} and \ref{rem:diff}). They allow the regularity of $F(u)$ to be any number below $p-1$ in the composition estimate, and below $p-2$ in the difference estimate. When $u$ is also controlled in $L^\infty$ or in a Lebesgue space, the classical theory gives estimates of $F(u)$ in Sobolev and Besov spaces of any order below $p$, see Ginibre, Ozawa and Velo \cite{GOV94} and \cite[Chapter 5]{RS96}. For the heat equation with this nonlinearity, the optimal Sobolev threshold for well-posedness was recently identified by Pineau and Taylor \cite{PT25}. These results do not apply directly in our setting. Indeed, let $0<s<\frac nr$ and $\frac1q=\frac1r-\frac sn$. By Proposition \ref{prop:weight}, the function $|x|^{-n/q}$ belongs to $\B sr\infty$, while it lies neither in $L^\infty$ nor in any space $L^\rho$ with $1\le\rho<\infty$. We do not know whether the thresholds $p-1$ and $p-2$ can be raised to $p$ and $p-1$ in the present setting. Such an improvement would replace $p-2$ by $p-1$ in \eqref{eq:H3}.

Finally, we compare the range of Theorem \ref{thm:global} with that of \cite[Theorem 3.10]{CF25}, which deals with $\gamma=0$ in the weak-Besov space $\dot B^{s,\infty}_{r,\infty}$.

\begin{remark}\label{rem:dimensions}
In our notation and with $\gamma=0$, the hypotheses of that theorem include $r<\frac{2(n+1)}{n-1}$, that is $2\beta_r<1$, and $s_r<1+\eta_r$. Its proof relies on Lemma 3.9 there, which also requires $p\ge2$. These conditions force
\[
(n-5)(p-1)<4 .
\]
Indeed, put $x:=\frac12-\frac1r$, so that $0<x<\frac1{n+1}$. By \eqref{eq:beta-eta} and \eqref{eq:indices} with $\gamma=0$, the condition $s_r<1+\eta_r$ is equivalent to $f(x)<0$, where
\[
f(x):=(p-1)\Bigl(\frac n2-1-(2n-1)x\Bigr)-1-(n-1)x .
\]
The function $f$ is decreasing, so $f(\frac1{n+1})<f(x)<0$. Since $f(\frac1{n+1})=\frac{n}{n+1}\bigl(\frac{(p-1)(n-5)}2-2\bigr)$, the claim follows. Together with $p\ge2$, this shows that the previous result applies only in dimensions $3\le n\le8$, and only to powers $p<1+\frac4{n-5}$ when $n=6,7,8$.

In contrast, for every $n\ge2$ the hypotheses of Theorem \ref{thm:global} hold for all sufficiently large $p$. Indeed, let $0<\gamma<\frac12$ and $r:=\frac{2(n+1)}n$. As in Remark \ref{rem:hypotheses} (c), $2\beta_r=\frac12$ and $\eta_r=\frac{n-1}{2(n+1)}$. Hence \eqref{eq:H2} holds, and \eqref{eq:H1} holds for $p>\frac{2(n+1)}{n-1}$. By \eqref{eq:sigma-formula},
\[
\lim_{p\to\infty}\sigma_r=\frac{n-2}2+\gamma+\frac1{2(n+1)}>0,
\]
while $p-2\to\infty$. Moreover, $\frac nr=\frac{n^2}{2(n+1)}\ge\frac23>\gamma$. Therefore \eqref{eq:H3} holds for every sufficiently large $p$.

\end{remark}

The paper is organized as follows. Section \ref{sec:linear} contains the linear estimates, namely the dispersive estimates for the wave group in homogeneous Besov spaces and the estimates for the singular weight, which rely on the product law in Besov space. Section \ref{sec:composition} is devoted to the composition and difference estimates for the fractional power and to their application to the weighted nonlinearity. In Section \ref{sec:wellposedness} we give a precise meaning to the time integrals, and we prove Theorem \ref{thm:global} and Corollary \ref{cor:data} by a contraction argument in the time-weighted space. Finally,  Section \ref{sec:asymptotic} contains the results on scattering, decay, asymptotic stability and asymptotic self-similarity. 

\section{Linear estimates}\label{sec:linear}

\subsection{Dispersive estimates for the wave group}

We start with the basic estimates for the wave group in homogeneous Besov spaces. The corresponding estimate in homogeneous Besov--Lorentz spaces is obtained in \cite[Lemma 3.1]{CF25}.

\begin{lemma}\label{lem:group}
Let $n\ge2$, $2\le r<\infty$, $s\in\R$, $1\le\sigma\le\infty$ and $0\le\theta\le1$. There is a constant $C>0$, depending on $n,r,\theta$, such that
\begin{align}
\|K(t)\psi\|_{\B sr\sigma}&\le C|t|^{-\theta\eta_r}\|\psi\|_{\B{s-1+(1-\theta)n(1-\frac2r)+2\theta\beta_r}{r'}{\sigma}},\label{eq:group-K}\\
\|\dK(t)\psi\|_{\B sr\sigma}&\le C|t|^{-\theta\eta_r}\|\psi\|_{\B{s+(1-\theta)n(1-\frac2r)+2\theta\beta_r}{r'}{\sigma}},\label{eq:group-dK}
\end{align}
for every $t\ne0$ and every $\psi$ in the space on the right-hand side. More precisely, for $j\in\Z$ and $t\ne0$,
\begin{equation}\label{eq:group-block}
\|\Dj jK(t)\psi\|_{L^r}\le C|t|^{-\theta\eta_r}2^{j\{-1+(1-\theta)n(1-\frac2r)+2\theta\beta_r\}}\sum_{|k-j|\le2}\|\Dj k\psi\|_{L^{r'}},
\end{equation}
and the same bound holds for $\dK(t)$ after removing $-1$ from the exponent. In particular, the choice $\theta=1$ gives
\begin{equation}\label{eq:group}
\|K(t)\psi\|_{\B sr\sigma}\le C|t|^{-\eta_r}\|\psi\|_{\B{s-1+2\beta_r}{r'}{\sigma}} .
\end{equation}
\end{lemma}

\begin{proof}
Let $U(t):=e^{it|D|}$. Then
\[
K(t)=\frac{U(t)-U(-t)}{2i|D|},\qquad \dK(t)=\frac{U(t)+U(-t)}2 .
\]
The stationary phase method yields the classical dispersive bound
\begin{equation}\label{eq:disp-1}
\|\Dj jU(t)f\|_{L^\infty}\le C\min\bigl\{2^{nj},|t|^{-\frac{n-1}2}2^{j\frac{n+1}2}\bigr\}\|f\|_{L^1}
\end{equation}
for every $j\in\Z$, $t\ne0$ and $f\in L^1(\Rn)$, with $C$ independent of $j$ and $t$ (see, e.g., \cite[Chapter 8]{BCD11} and \cite{GV89}). By the Plancherel theorem we also have $\|\Dj jU(t)f\|_{L^2}\le\|f\|_{L^2}$. The Riesz--Thorin theorem, applied with $\frac1r=\frac{1-2/r}\infty+\frac{2/r}2$, therefore gives
\[
\|\Dj jU(t)f\|_{L^r}\le C\min\bigl\{2^{nj(1-\frac2r)},|t|^{-\eta_r}2^{2j\beta_r}\bigr\}\|f\|_{L^{r'}},
\]
because $\frac{n-1}2(1-\frac2r)=\eta_r$ and $\frac{n+1}2(1-\frac2r)=2\beta_r$ by \eqref{eq:beta-eta}. We only need an upper bound, so we may use $\min\{A,B\}\le A^{1-\theta}B^\theta$. This gives
\begin{equation}\label{eq:disp-2}
\|\Dj jU(t)f\|_{L^r}\le C|t|^{-\theta\eta_r}2^{j\{(1-\theta)n(1-\frac2r)+2\theta\beta_r\}}\|f\|_{L^{r'}} .
\end{equation}
Next, $\Dj j=\Dj j\tDj j$ and $|D|^{-1}$ commutes with $\Dj j$, so that
\[
\Dj jK(t)\psi=\frac1{2i}\bigl(\Dj jU(t)-\Dj jU(-t)\bigr)\bigl(|D|^{-1}\tDj j\psi\bigr).
\]
Let $\widetilde\phi:=\sum_{|l|\le1}\phi(2^{-l}\cdot)$, which equals one on $\supp\phi$, and $\tilde g:=\F^{-1}[|\xi|^{-1}\widetilde\phi]\in\Sch(\Rn)$. Then
\[
|D|^{-1}\Dj k\psi=\bigl(2^{k(n-1)}\tilde g(2^k\cdot)\bigr)*\Dj k\psi,\qquad\bigl\|2^{k(n-1)}\tilde g(2^k\cdot)\bigr\|_{L^1}=2^{-k}\|\tilde g\|_{L^1},
\]
and the Young inequality gives $\||D|^{-1}\Dj k\psi\|_{L^{r'}}\le C2^{-k}\|\Dj k\psi\|_{L^{r'}}$. We apply \eqref{eq:disp-2} to $f=|D|^{-1}\tDj j\psi$, which is the sum of $|D|^{-1}\Dj k\psi$ over $|k-j|\le2$, and obtain \eqref{eq:group-block}. After multiplication by $2^{js}$ this becomes
\[
2^{js}\|\Dj jK(t)\psi\|_{L^r}\le C|t|^{-\theta\eta_r}\sum_{|k-j|\le2}2^{k\{s-1+(1-\theta)n(1-\frac2r)+2\theta\beta_r\}}\|\Dj k\psi\|_{L^{r'}},
\]
and the $\ell^\sigma$ norm in $j$ gives \eqref{eq:group-K}. The proof of \eqref{eq:group-dK} is identical, except that the factor $|D|^{-1}$ is absent.
\end{proof}

\begin{remark}\label{rem:group}
The fixed point argument uses only \eqref{eq:group}, through Lemma \ref{lem:group-potential} below. Estimates \eqref{eq:group-K} and \eqref{eq:group-dK} with $\theta<1$ serve to describe the class of initial data (Corollary \ref{cor:data}). Related dispersive estimates in homogeneous Besov and Sobolev spaces appear in \cite{GV89,Pec00,RY02,FPL20}.
\end{remark}

\subsection{Linear estimates involving the singular weight}

We recall the following fact.

\begin{proposition}\label{prop:weight}
Let $0<b<n$ and $\Omega\in C^\infty(\mathbb S^{n-1})$, and put $f(x):=\Omega(x/|x|)|x|^{-b}$. For every $m\in[1,\infty]$, the function $f$ belongs to $\B{\frac nm-b}m\infty(\Rn)$. In particular, $h_\gamma=\abs{\cdot}^{-\gamma}\in\B{\frac nm-\gamma}m\infty(\Rn)$ for $0<\gamma<n$ (see \cite[p.~65]{BCD11}).
\end{proposition}

The product law of Abidi \cite[Proposition 2.2]{Abi07}, recalled in Proposition \ref{prop:Abidi} of the Appendix, leads to the following estimate. In this law one factor is only required to belong to a homogeneous Besov space with third index $\infty$.

\begin{lemma}\label{lem:potential}
Let $h_\gamma=\abs{\cdot}^{-\gamma}$ with $0<\gamma<n$. Assume that $1\le m,\sigma\le\infty$ and $\varsigma\in\R$ satisfy
\begin{equation}\label{eq:potential-cond}
\gamma-\frac nm-n\min\Bigl\{0,1-\frac2m\Bigr\}<\varsigma<\frac nm .
\end{equation}
Then there is a constant $C>0$, depending on $n,\varsigma,\sigma,m$ and $\gamma$, such that
\begin{equation}\label{eq:potential}
\|h_\gamma f\|_{\B{\varsigma-\gamma}m\sigma}\le C\|h_\gamma\|_{\B{\frac nm-\gamma}m\infty}\|f\|_{\B \varsigma m\sigma}.
\end{equation}
In particular, if $2\le r\le\infty$ and $\gamma-\frac nr<\varsigma<\frac n{r'}$, then $\|h_\gamma f\|_{\B{\varsigma-\gamma}{r'}\sigma}\le C\|f\|_{\B\varsigma{r'}\sigma}$.
\end{lemma}

\begin{proof}
Proposition \ref{prop:weight} gives $h_\gamma\in\B{\frac nm-\gamma}m\infty$. We apply Proposition \ref{prop:Abidi} with $s_1:=\varsigma$ and $s_2:=\frac nm-\gamma$, so that $s_1+s_2-\frac nm=\varsigma-\gamma$. The requirements $s_1<\frac nm$ and $s_1+s_2+n\min\{0,1-\frac2m\}>0$ are exactly \eqref{eq:potential-cond}, and $s_2<\frac nm$ holds because $\gamma>0$. For the last assertion, observe that $m=r'\le2$. Then $n\min\{0,1-\frac2{r'}\}=n(\frac2r-1)$, and the left-hand side of \eqref{eq:potential-cond} equals $\gamma-\frac n{r'}-\frac{2n}r+n=\gamma-\frac nr$.
\end{proof}

We next prove dispersive estimates for the inhomogeneous term that contains the singular weight. We first treat the time $t=1$, and then we use a scaling argument.

\begin{lemma}\label{lem:group-potential}
Let $h_\gamma=\abs{\cdot}^{-\gamma}$ with $0<\gamma<n$ as in Lemma \ref{lem:potential}. Let $2\le r<\infty$, $1\le\sigma\le\infty$ and $0\le\theta\le1$, and let $\varsigma\in\R$ satisfy
\begin{equation}\label{eq:group-potential-cond}
\gamma-\frac nr<\varsigma<\frac n{r'} .
\end{equation}
Put
\[
\tilde s_\theta:=\varsigma-\gamma+1-(1-\theta)n\Bigl(1-\frac2r\Bigr)-2\theta\beta_r .
\]
Then there is a constant $C>0$, depending on $n,r,\theta,\varsigma,\sigma$ and $\gamma$, such that
\begin{equation}\label{eq:group-potential}
\|K(t)(h_\gamma u)\|_{\B{\tilde s_\theta}r\sigma}+\|\dK(t)(h_\gamma u)\|_{\B{\tilde s_\theta-1}r\sigma}\le C|t|^{-\theta\eta_r}\|u\|_{\B\varsigma{r'}\sigma}
\end{equation}
for all $t\ne0$ and $u\in\B\varsigma{r'}\sigma(\Rn)$. In particular, the choice $\theta=1$ gives
\begin{equation}\label{eq:group-potential-1}
\|K(t)(h_\gamma u)\|_{\B{\varsigma-\gamma+1-2\beta_r}r\sigma}\le C|t|^{-\eta_r}\|u\|_{\B\varsigma{r'}\sigma}.
\end{equation}
\end{lemma}

\begin{proof}
We first take $t=1$. Since $\tilde s_\theta-1+(1-\theta)n(1-\frac2r)+2\theta\beta_r=\varsigma-\gamma$, Lemma \ref{lem:group} with $t=1$ and $s=\tilde s_\theta$ gives
\[
\|K(1)\psi\|_{\B{\tilde s_\theta}r\sigma}+\|\dK(1)\psi\|_{\B{\tilde s_\theta-1}r\sigma}\le C\|\psi\|_{\B{\varsigma-\gamma}{r'}\sigma}
\]
for all $\psi\in\B{\varsigma-\gamma}{r'}\sigma$. By Lemma \ref{lem:potential} with $m=r'$, there is a constant $C$ depending on $n,\varsigma,\sigma,r$ and $\gamma$ such that
\[
\|h_\gamma u\|_{\B{\varsigma-\gamma}{r'}\sigma}\le C\|h_\gamma\|_{\B{\frac n{r'}-\gamma}{r'}\infty}\|u\|_{\B\varsigma{r'}\sigma},
\]
provided that \eqref{eq:group-potential-cond} holds. This proves \eqref{eq:group-potential} for $t=1$.

We now remove the restriction $t=1$ by scaling. Since $K(-t)=-K(t)$ and $\dK(-t)=\dK(t)$, we may assume $t>0$. For $\lambda>0$ we write $\delta_\lambda f:=f(\lambda\,\cdot)$. A computation on the Fourier side gives $K(t)\delta_\lambda=\lambda^{-1}\delta_\lambda K(\lambda t)$ and $\dK(t)\delta_\lambda=\delta_\lambda\dK(\lambda t)$. Moreover, $\delta_{1/t}(h_\gamma\,\delta_tu)=t^\gamma h_\gamma u$. With $\lambda=1/t$ we obtain
\[
K(t)(h_\gamma u)=t^{1-\gamma}\delta_{1/t}\bigl[K(1)(h_\gamma\,\delta_tu)\bigr],\qquad
\dK(t)(h_\gamma u)=t^{-\gamma}\delta_{1/t}\bigl[\dK(1)(h_\gamma\,\delta_tu)\bigr].
\]
We use the dilation bound $\|\delta_\lambda f\|_{\B s{\rho}\sigma}\le C\lambda^{s-\frac n\rho}\|f\|_{\B s\rho\sigma}$, with $C$ independent of $\lambda>0$ (see \cite[Chapter 2]{BCD11}), and the case $t=1$. This gives
\[
\|K(t)(h_\gamma u)\|_{\B{\tilde s_\theta}r\sigma}\le Ct^{1-\gamma-\tilde s_\theta+\frac nr}\|K(1)(h_\gamma\,\delta_tu)\|_{\B{\tilde s_\theta}r\sigma}
\le Ct^{1-\gamma-\tilde s_\theta+\frac nr+\varsigma-\frac n{r'}}\|u\|_{\B\varsigma{r'}\sigma}.
\]
By the definition of $\tilde s_\theta$, the exponent of $t$ equals
\[
(1-\theta)n\Bigl(1-\frac2r\Bigr)+2\theta\beta_r-n\Bigl(1-\frac2r\Bigr)=-\theta\Bigl(n\Bigl(1-\frac2r\Bigr)-2\beta_r\Bigr)=-\theta\eta_r,
\]
because $n(1-\frac2r)-2\beta_r=(n-1)(\frac12-\frac1r)=\eta_r$ by \eqref{eq:beta-eta}. The same computation applies to $\dK(t)$, where the factor $t^{-\gamma}$ and the regularity $\tilde s_\theta-1$ give the same power of $t$. This proves \eqref{eq:group-potential}, and \eqref{eq:group-potential-1} is the case $\theta=1$.
\end{proof}

\section{Estimates of the fractional power nonlinearity in Besov spaces}\label{sec:composition}

For $0<\rho\le\infty$ and a compact set $\Omega\subset\Rn$ we write
\[
L^\rho_\Omega(\Rn):=\{f\in\Sch'(\Rn):\ \supp\widehat f\subset\Omega,\ \|f\|_{L^\rho}<\infty\}.
\]
We recall the Nikol'skii inequality, a version of the Bernstein inequality that remains valid for exponents below one.

\begin{proposition}[{\cite[p.~22]{Tri83}}]\label{prop:Bernstein}
Let $0<\rho_1\le\rho_2\le\infty$ and let $\alpha$ be a multi-index. Then there is a constant $C>0$ such that
\[
\|\partial^\alpha f\|_{L^{\rho_2}}\le CR^{|\alpha|+n(\frac1{\rho_1}-\frac1{\rho_2})}\|f\|_{L^{\rho_1}}
\]
for all $R>0$ and all $f\in L^{\rho_1}_{B(0,R)}(\Rn)$, where $B(0,R):=\{|\xi|\le R\}$.
\end{proposition}

(The case $R=1$ is \cite[p.~22]{Tri83}, and the general case follows by scaling.) Besides Proposition \ref{prop:Bernstein}, we use its reverse form on annuli \cite[Lemma 2.1]{BCD11}, written for Littlewood--Paley blocks. Let $N\ge0$ be an integer. Since $|\xi|^{2N}=\sum_{|\alpha|=N}\frac{N!}{\alpha!}\xi^{2\alpha}$, we have
\[
\phi_l(\xi)=2^{-lN}\sum_{|\alpha|=N}(i\xi)^\alpha\,\widehat{\chi_\alpha}(2^{-l}\xi),\qquad
\widehat{\chi_\alpha}(\eta):=\frac{N!}{\alpha!}(-i)^N\frac{\eta^\alpha\phi(\eta)}{|\eta|^{2N}} .
\]
Each $\chi_\alpha$ belongs to $\Sch(\Rn)$, because $\phi$ vanishes near the origin, and $\int\chi_\alpha\,dx=\widehat{\chi_\alpha}(0)=0$. Let $1\le\vartheta\le\infty$, and let $g\in\Sch'(\Rn)$ be such that its derivatives of order $N$ belong to $L^\vartheta(\Rn)$. Then
\[
\Dj lg=2^{-lN}\sum_{|\alpha|=N}\chi_{\alpha,l}*\partial^\alpha g,\qquad\chi_{\alpha,l}:=2^{ln}\chi_\alpha(2^l\cdot),
\]
and the Young inequality gives the reverse Bernstein inequality
\begin{equation}\label{eq:reverse}
\|\Dj lg\|_{L^\vartheta}\le C_N2^{-lN}\sup_{|\alpha|=N}\|\partial^\alpha g\|_{L^\vartheta}\qquad\text{for all }l\in\Z,
\end{equation}
where $C_N:=\sum_{|\alpha|=N}\|\chi_\alpha\|_{L^1}$. For $N=0$ this is the boundedness of $\Dj l$ on $L^\vartheta(\Rn)$. For $y\in\Rn$ we write $\tau_yh:=h(\cdot-y)$. Since $\int\chi_\alpha\,dx=0$, we also have $\chi_{\alpha,l}*h=\int_{\Rn}\chi_{\alpha,l}(y)(\tau_yh-h)\,dy$, and the Minkowski inequality refines \eqref{eq:reverse} as follows. For $0\le\beta\le1$, $A\ge0$ and $k\in\Z$,
\begin{equation}\label{eq:difference-char}
\sup_{|\alpha|=N}\|\tau_y\partial^\alpha g-\partial^\alpha g\|_{L^\vartheta}\le A\,2^{kN}(2^k|y|)^\beta\ \ \forall\,y\in\Rn
\ \Longrightarrow\
\|\Dj lg\|_{L^\vartheta}\le C_{N,\beta}A\,2^{(k-l)(N+\beta)}\ \ \forall\,l\in\Z,
\end{equation}
where $C_{N,\beta}:=\sum_{|\alpha|=N}\int_{\Rn}|\chi_\alpha(z)||z|^\beta dz$. For $N=0$ this is the zero-mean argument in the proof of \cite[Theorem 2.36]{BCD11}, which characterizes $\B s\vartheta\sigma$, $0<s<1$, by first differences. Thus $N$ derivatives produce the decay $2^{(k-l)N}$ through the reverse Bernstein inequality, and one first difference of the derivatives of order $N$ produces the remaining decay $2^{(k-l)\beta}$. We shall also use the following consequence of Proposition \ref{prop:Bernstein}. If $1\le\rho\le\infty$, $f\in L^\rho_{B(0,R)}$ and $y\in\Rn$, then $\|\tau_yf-f\|_{L^\rho}\le2\|f\|_{L^\rho}$ and $\|\tau_yf-f\|_{L^\rho}\le|y|\|\nabla f\|_{L^\rho}\le CR|y|\|f\|_{L^\rho}$, and therefore
\begin{equation}\label{eq:translation}
\|\tau_yf-f\|_{L^\rho}\le C\min\{1,R|y|\}\|f\|_{L^\rho}.
\end{equation}

Next we collect elementary facts on $F(z)=|z|^{p-1}z$. Let $m\ge1$ be an integer with $m<p$. Then $F\in C^m(\R)$, and
\[
F^{(m)}(z)=\begin{cases}c_m|z|^{p-m}&\text{if $m$ is odd},\\ c_m|z|^{p-m-1}z&\text{if $m$ is even},\end{cases}\qquad c_m:=p(p-1)\cdots(p-m+1).
\]
This follows by induction on $m$, since for $\mu>1$ the functions $|z|^\mu$ and $|z|^{\mu-1}z$ are continuously differentiable, with derivatives $\mu|z|^{\mu-2}z$ and $\mu|z|^{\mu-1}$. In particular, $|F^{(m)}(\omega)|\le|c_m||\omega|^{p-m}$. Moreover, if $0\le\beta\le\min\{1,p-m\}$, then
\begin{equation}\label{eq:F1}
|F^{(m)}(\omega_1)-F^{(m)}(\omega_2)|\le C(|\omega_1|+|\omega_2|)^{p-m-\beta}|\omega_1-\omega_2|^\beta,\qquad\omega_1,\omega_2\in\R.
\end{equation}
For $p=2$ we put $F''(z):=2\operatorname{sgn}z$, with $F''(0):=0$. Then, for every $p\ge2$, the function $F'$ is locally Lipschitz and
\begin{equation}\label{eq:F3}
\begin{aligned}
&F'(\omega_1)-F'(\omega_2)=(\omega_1-\omega_2)\int_0^1F''\bigl(\omega_2+\lambda(\omega_1-\omega_2)\bigr)d\lambda,\\
&\Bigl|\int_0^1F''\bigl(\omega_2+\lambda(\omega_1-\omega_2)\bigr)d\lambda\Bigr|\le C(|\omega_1|+|\omega_2|)^{p-2},\qquad\omega_1,\omega_2\in\R.
\end{aligned}
\end{equation}
To check \eqref{eq:F1}, let $g(\omega)=|\omega|^{\mu}$ or $g(\omega)=|\omega|^{\mu-1}\omega$ with $\mu>0$. For $\mu\ge1$ the mean value theorem gives $|g(\omega_1)-g(\omega_2)|\le\mu(|\omega_1|+|\omega_2|)^{\mu-1}|\omega_1-\omega_2|$, while for $0<\mu<1$ the function $g$ is Hölder continuous of order $\mu$, namely $|g(\omega_1)-g(\omega_2)|\le2^{1-\mu}|\omega_1-\omega_2|^\mu$. We combine these bounds with $|g(\omega_1)-g(\omega_2)|\le2(|\omega_1|+|\omega_2|)^\mu$ and with $\min\{A_1,A_2\}\le A_1^{1-\vartheta}A_2^\vartheta$ ($0\le\vartheta\le1$). This yields \eqref{eq:F1} with $\mu=p-m$.

Using the telescopic representation, gives
\begin{equation}\label{eq:telescopic}
\begin{aligned}
F(u)&=\lim_{N\to\infty}\bigl(F(\Sj Nu)-F(\Sj{-N}u)\bigr)
=\lim_{N\to\infty}\sum_{k=-N+1}^{N}\bigl(F(\Sj ku)-F(\Sj{k-1}u)\bigr)\\
&=\sum_{k\in\Z}\int_0^1F'\bigl(\Sj{k-1}u+t\Dj ku\bigr)dt\ \Dj ku=\sum_{k\in\Z}M_k(u)\,\Dj ku,
\end{aligned}
\end{equation}
where
\begin{equation}\label{eq:Mk}
M_k(u):=\int_0^1F'\bigl(z^t_k(u)\bigr)dt,\qquad z^t_k(u):=\Sj{k-1}u+t\Dj ku,\qquad k\in\Z,\ 0\le t\le1.
\end{equation}

\subsection{High frequencies of the coefficients}

 The next lemma adds the decay of their high frequencies. It is a new ingredient in the proofs of Lemmas \ref{lem:composition} and \ref{lem:difference} below. The idea comes from the composition theorem for smooth functions \cite[Theorem 2.61]{BCD11}. There, the coefficients $m_j=\int_0^1f'(\Sj ju+t\Dj ju)\,dt$ of the first linearization satisfy $\|\partial^\alpha m_j\|_{L^\infty}\le C2^{j|\alpha|}$ for every multi-index $\alpha$ \cite[Lemma 2.63]{BCD11}, and these bounds are turned into Besov regularity by \cite[Lemma 2.64]{BCD11}. Lemma \ref{lem:gain} is a version of \cite[Lemma 2.63]{BCD11} for the nonlinearity $|u|^{p-1}u$, which has only finite smoothness. It differs in three aspects. The coefficients are measured in $L^\theta$ instead of $L^\infty$. The derivatives are converted into decay of the high frequencies by the reverse Bernstein inequality \eqref{eq:reverse}. The fractional part of the gain comes from the Hölder continuity of the last derivative of $F$, through \eqref{eq:difference-char}. Hölder estimates of this kind for the derivatives of $|u|^{\mu}u$ are classical in the theory of Nemytskij operators, see \cite[Chapter 5]{RS96}.

\begin{lemma}\label{lem:gain}
Let $n\ge1$, $F(u):=|u|^{p-1}u$ with $\R\ni p>1$, and let $s\in\R$ and $r,\theta\in[1,\infty]$ satisfy
\begin{equation}\label{eq:gain-cond}
\kappa:=\frac nr-s>0,\qquad r\le(p-1)\theta,\qquad\frac n\theta<(p-1)\kappa .
\end{equation}
Let $M_k$ be given by \eqref{eq:Mk}, and write $\|\cdot\|:=\|\cdot\|_{\B sr\infty}$.
\begin{enumerate}[(i)]
\item Let $0\le\nu\le p-1$. There is a constant $C>0$ such that, for all $u\in\B sr\infty$ and $k,l\in\Z$,
\begin{equation}\label{eq:Mk-gain}
\|M_k(u)\|_{L^\theta}+2^{(l-k)\nu}\|\Dj lM_k(u)\|_{L^\theta}\le C2^{k((p-1)\kappa-\frac n\theta)}\|u\|^{p-1}.
\end{equation}
\item Let $p\ge2$ and $0\le\nu\le p-2$. There is a constant $C>0$ such that, for all $u,v\in\B sr\infty$ and $k,l\in\Z$, the function $\widetilde M_k:=M_k(u)-M_k(v)$ satisfies
\begin{equation}\label{eq:tildeM-gain}
\|\widetilde M_k\|_{L^\theta}+2^{(l-k)\nu}\|\Dj l\widetilde M_k\|_{L^\theta}\le C2^{k((p-1)\kappa-\frac n\theta)}\bigl(\|u\|+\|v\|\bigr)^{p-2}\|u-v\|.
\end{equation}
\end{enumerate}
\end{lemma}

\begin{proof}
First, we try to obtain bounds in Lebesgue spaces. Let $\varrho:=(p-1)\theta$. By \eqref{eq:gain-cond}, $r\le\varrho\le\infty$. For $0\le t\le1$, Proposition \ref{prop:Bernstein} gives
\begin{equation}\label{eq:a-bound}
\begin{aligned}
\|z^t_k(u)\|_{L^\varrho}&\le\sum_{k'\le k}\|\Dj{k'}u\|_{L^\varrho}\le C\sum_{k'\le k}2^{k'n(\frac1r-\frac1\varrho)}\|\Dj{k'}u\|_{L^r}
=C\sum_{k'\le k}2^{k'(\kappa-\frac n\varrho)}2^{k's}\|\Dj{k'}u\|_{L^r}\\
&\le C2^{k(\kappa-\frac n\varrho)}\|u\|,
\end{aligned}
\end{equation}
because $\kappa-\frac n\varrho=\frac1{p-1}\bigl((p-1)\kappa-\frac n\theta\bigr)>0$ by \eqref{eq:gain-cond}. More generally, for every multi-index $\alpha$ and every $y\in\Rn$,
\begin{equation}\label{eq:a-translation}
\begin{aligned}
\|\partial^\alpha z^t_k(u)\|_{L^\varrho}&\le C2^{k(|\alpha|+\kappa-\frac n\varrho)}\|u\|,\\
\|\tau_y\partial^\alpha z^t_k(u)-\partial^\alpha z^t_k(u)\|_{L^\varrho}&\le C\min\{1,2^k|y|\}\,2^{k(|\alpha|+\kappa-\frac n\varrho)}\|u\|.
\end{aligned}
\end{equation}
Indeed, Proposition \ref{prop:Bernstein} gives $\|\partial^\alpha\Dj{k'}u\|_{L^\varrho}\le C2^{k'(|\alpha|+\kappa-\frac n\varrho)}\|u\|$ as in \eqref{eq:a-bound}, and \eqref{eq:translation} adds the factor $\min\{1,2^{k'}|y|\}$ to this bound. We sum over $k'\le k$ and use $\sum_{k'\le k}\min\{1,2^{k'}|y|\}2^{k'b}\le C\min\{1,2^k|y|\}2^{kb}$, which holds for $b>0$. The same computation with $\varrho=\infty$ shows that $z^t_k(u)$ is a smooth function whose derivatives of all orders are bounded. Since $|F'(z)|=p|z|^{p-1}$ and $\frac{p-1}\varrho=\frac1\theta$, we deduce from \eqref{eq:a-bound} that
\begin{equation}\label{eq:Mk-bound}
\|M_k(u)\|_{L^\theta}\le\int_0^1\|F'(z^t_k(u))\|_{L^\theta}dt\le p\sup_{0\le t\le1}\|z^t_k(u)\|_{L^\varrho}^{p-1}\le C2^{k((p-1)\kappa-\frac n\theta)}\|u\|^{p-1}.
\end{equation}

Now, we deal with high frequencies of $M_k(u)$. $\nu=0$, then \eqref{eq:Mk-gain} follows from \eqref{eq:Mk-bound}, since $\Dj l$ is bounded on $L^\theta(\Rn)$. Let $\nu>0$, and write $\nu=N+\beta$ with $N:=\lceil\nu\rceil-1$ and $\beta\in(0,1]$. Since $\nu\le p-1$, we have $N+1<p$, and $\beta\le\min\{1,p-1-j\}$ for $0\le j\le N$. In particular, $F\in C^{N+1}(\R)$. Write $z=z^t_k(u)$ and $\omega_k:=2^{k((p-1)\kappa-\frac n\theta)}\|u\|^{p-1}$. Since $z$ is smooth with bounded derivatives, $F'(z)$ is $N$ times continuously differentiable.

Fix a multi-index $\alpha$ with $|\alpha|=N$. By the chain rule (the Faà di Bruno formula), as in the proof of \cite[Lemma 2.63]{BCD11}, $\partial^\alpha(F'(z))$ is a linear combination, with coefficients depending only on $\alpha$, of the functions
\[
P:=F^{(1+j)}(z)\,\partial^{\alpha_1}z\cdots\partial^{\alpha_j}z,\qquad 0\le j\le N,\quad \alpha_1+\dots+\alpha_j=\alpha,\quad|\alpha_i|\ge1,
\]
where $j=0$, that is $P=F'(z)$, occurs only for $N=0$. For such a term,
\[
\tau_yP-P=\bigl(F^{(1+j)}(\tau_yz)-F^{(1+j)}(z)\bigr)\prod_{i=1}^j\tau_y\partial^{\alpha_i}z+F^{(1+j)}(z)\Bigl(\prod_{i=1}^j\tau_y\partial^{\alpha_i}z-\prod_{i=1}^j\partial^{\alpha_i}z\Bigr).
\]
The last difference is a sum of $j$ products. Each of them contains one factor $\tau_y\partial^{\alpha_i}z-\partial^{\alpha_i}z$, and its other factors are of the form $\partial^{\alpha_{i'}}z$ or $\tau_y\partial^{\alpha_{i'}}z$. The first term is controlled by \eqref{eq:F1} with $m=1+j$, and the second one by $|F^{(1+j)}(z)|\le C|z|^{p-1-j}$. In every resulting product the exponents of the factors add up to $(p-1-j-\beta)+\beta+j=p-1$ or to $(p-1-j)+j=p-1$, and $\frac{p-1}\varrho=\frac1\theta$. The generalized Hölder inequality (see \cite[p.~10]{Gra08}), \eqref{eq:a-bound}, \eqref{eq:a-translation}, $|\alpha_1|+\dots+|\alpha_j|=N$ and $\min\{1,x\}\le x^\beta$ therefore give
\[
\|\tau_yP-P\|_{L^\theta}\le C\bigl(\min\{1,2^k|y|\}^\beta+\min\{1,2^k|y|\}\bigr)2^{kN}\omega_k\le C2^{kN}(2^k|y|)^\beta\omega_k .
\]
The same argument without increments shows that $\|P\|_{L^\theta}\le C2^{kN}\omega_k$. Moreover, $\partial^\alpha M_k(u)$ is the integral over $t\in[0,1]$ of $\partial^\alpha\bigl(F'(z^t_k(u))\bigr)$, and $\theta\ge1$. Hence $\partial^\alpha M_k(u)\in L^\theta(\Rn)$ and
\[
\|\tau_y\partial^\alpha M_k(u)-\partial^\alpha M_k(u)\|_{L^\theta}\le C2^{kN}(2^k|y|)^\beta\omega_k,\qquad|\alpha|=N,\ y\in\Rn .
\]
Therefore \eqref{eq:difference-char}, applied with $g=M_k(u)$, $\vartheta=\theta$ and $A=C\omega_k$, gives $\|\Dj lM_k(u)\|_{L^\theta}\le C2^{(k-l)\nu}\omega_k$ for all $l\in\Z$. Together with \eqref{eq:Mk-bound}, this proves (i).

Now, we compute bounds for $\widetilde M_k$ in Lebesgue spaces. Let $p\ge2$, put $w:=u-v$, $z_1:=z^t_k(u)$ and $z_2:=z^t_k(v)$. By linearity, $z_1-z_2=z^t_k(w)=\Sj{k-1}w+t\Dj kw$. From \eqref{eq:F3} we get
\begin{equation}\label{eq:tildeM-repr}
\widetilde M_k=\int_0^1(z_1-z_2)\,G\,dt,\qquad G:=\int_0^1F''\bigl(z_2+\lambda(z_1-z_2)\bigr)d\lambda,\qquad |G|\le C(|z_1|+|z_2|)^{p-2}.
\end{equation}
The Hölder inequality with $\frac1\varrho+\frac{p-2}\varrho=\frac1\theta$ and \eqref{eq:a-bound}, applied to $z_1$, $z_2$ and $z_1-z_2$, then give
\begin{equation}\label{eq:tildeM-bound}
\|\widetilde M_k\|_{L^\theta}\le C\sup_{0\le t\le1}\|z_1-z_2\|_{L^\varrho}\bigl(\|z_1\|_{L^\varrho}+\|z_2\|_{L^\varrho}\bigr)^{p-2}
\le C2^{k((p-1)\kappa-\frac n\theta)}\bigl(\|u\|+\|v\|\bigr)^{p-2}\|w\|.
\end{equation}

Finally, we  obtain high frequencies of $\widetilde M_k$. In fact, if $\nu=0$, then \eqref{eq:tildeM-gain} follows from \eqref{eq:tildeM-bound}. Let $\nu>0$, and write $\nu=N+\beta$ with $N:=\lceil\nu\rceil-1$ and $\beta\in(0,1]$. Since $\nu\le p-2$, we have $N+2<p$, and $\beta\le\min\{1,p-2-j\}$ for $0\le j\le N$. In particular, $F\in C^{N+2}(\R)$. Put $d:=z_1-z_2$, $\varpi_\lambda:=z_2+\lambda d$ and $\widetilde\omega_k:=2^{k((p-1)\kappa-\frac n\theta)}(\|u\|+\|v\|)^{p-2}\|w\|$. Then $|\varpi_\lambda|\le|z_1|+|z_2|$, and $\partial^{\alpha}\varpi_\lambda=(1-\lambda)\partial^{\alpha}z_2+\lambda\partial^{\alpha}z_1$ for every multi-index $\alpha$. Moreover, \eqref{eq:a-bound} and \eqref{eq:a-translation} apply to $z_1$, $z_2$ and $d=z^t_k(w)$, with $u$, $v$ and $w$ respectively.

Fix $\alpha$, with $|\alpha|=N$. By the Leibniz rule and the chain rule, $\partial^\alpha(dG)$ is a linear combination of the functions
\[
Q:=\partial^{\alpha_0}d\int_0^1F^{(2+j)}(\varpi_\lambda)\,\partial^{\alpha_1}\varpi_\lambda\cdots\partial^{\alpha_j}\varpi_\lambda\,d\lambda,
\]
where $0\le j\le N$, $\alpha_0+\dots+\alpha_j=\alpha$ and $|\alpha_i|\ge1$ for $i\ge1$.
As for $P$ above, we write $\tau_yQ-Q$ as a sum of terms in which exactly one factor is replaced by its increment. The increment of $F^{(2+j)}(\varpi_\lambda)$ is controlled by \eqref{eq:F1} with $m=2+j$, together with $|\tau_y\varpi_\lambda-\varpi_\lambda|\le|\tau_yz_1-z_1|+|\tau_yz_2-z_2|$. The other factors are controlled by $|F^{(2+j)}(\varpi_\lambda)|\le C(|z_1|+|z_2|)^{p-2-j}$, by $|\partial^{\alpha_i}\varpi_\lambda|\le|\partial^{\alpha_i}z_1|+|\partial^{\alpha_i}z_2|$, and by the analogous bounds for the translates and their increments. In every resulting product the exponents of the factors add up to $1+(p-2-j)+j=p-1$. The Hölder inequality, \eqref{eq:a-bound}, \eqref{eq:a-translation}, $|\alpha_0|+\dots+|\alpha_j|=N$ and $\min\{1,x\}\le x^\beta$ give
\[
\|\tau_y\partial^\alpha(dG)-\partial^\alpha(dG)\|_{L^\theta}\le C2^{kN}(2^k|y|)^\beta\widetilde\omega_k .
\]
Integrating in $t$, we obtain the same bound for $\tau_y\partial^\alpha\widetilde M_k-\partial^\alpha\widetilde M_k$. Therefore \eqref{eq:difference-char} gives $\|\Dj l\widetilde M_k\|_{L^\theta}\le C2^{(k-l)\nu}\widetilde\omega_k$ for all $l\in\Z$. Together with \eqref{eq:tildeM-bound}, this proves (ii).
\end{proof}
The first bound in \eqref{eq:Mk-gain}, namely the estimate of $\|M_k(u)\|_{L^\theta}$, is \cite[(3.6)]{Chi19}, with $b=p-1$, $q=r$, $\theta_1=\theta$ and $d=n$ in the notation used there. Indeed, $b\{d(\frac1q-\frac1{b\theta_1})-s\}=(p-1)\kappa-\frac n\theta$, and the conditions $\frac1{b\theta_1}\le\frac1q$ and $s<d(\frac1q-\frac1{b\theta_1})$ required there are the last two conditions in \eqref{eq:gain-cond}. The first bound in \eqref{eq:tildeM-gain} is used in the proof of Lemma 3.7 there, through the mean value formula (3.14) there, without being written out. For $\nu=0$ the whole lemma therefore follows from that paper, since $\Dj l$ is bounded on $L^\theta(\Rn)$. The new part is the decay of $\Dj lM_k(u)$ and $\Dj l\widetilde M_k$ for $\nu>0$. The parabolic argument does not need it, because there the nonlinearity is measured in Besov spaces of negative regularity, and the term $III_j$ is estimated with $\|M_k(u)\|_{L^{\theta_2}}$ only, see (3.11) there. The size of the gain $\nu$ is limited only by the smoothness of $F$. The proof uses $N:=\lceil\nu\rceil-1$ derivatives of $F'(z^t_k(u))$, through \eqref{eq:reverse}, and the Hölder continuity of order $\nu-N$ of $F^{(1+N)}$, through \eqref{eq:difference-char}. If $\nu<p-1$ is an integer, one may apply \eqref{eq:reverse} with $N=\nu$ instead, and then the proof uses only Proposition \ref{prop:Bernstein}, its reverse form and the chain rule. If $p$ is an odd integer, then $F$ is a polynomial, and \eqref{eq:Mk-gain} holds for every $\nu\ge0$. Part (ii) requires $p\ge2$, because it uses \eqref{eq:F3}.

\subsection{Composition estimate}

The next lemma is a composition estimate for the fractional power in homogeneous Besov space, in the spirit of \cite[Lemma 3.2]{Chi19}. We state it in more generality than needed for Theorem \ref{thm:global}, since it may be useful elsewhere.

\begin{lemma}\label{lem:composition}
Let $n\ge1$, $F(u):=|u|^{p-1}u$ with $\R\ni p>1$, $1\le\sigma\le\infty$ and $0\le\nu\le p-1$. Assume that $s\in\R$ and $r,m,\theta_1,\theta_2\in[1,\infty]$ satisfy, for $i=1,2$,
\begin{equation}\label{eq:comp-structural}
\kappa:=\frac nr-s>0,\qquad r\le(p-1)\theta_i,\qquad \frac n{\theta_i}<(p-1)\kappa,\qquad \frac1m\le\frac1{r}+\frac1{\theta_i},
\end{equation}
and define
\begin{equation}\label{eq:s1}
s_1:=\frac nm-p\kappa=\frac nm-p\Bigl(\frac nr-s\Bigr).
\end{equation}
Assume moreover that
\begin{equation}\label{eq:comp-highhigh}
\begin{cases}
s_1+n\Bigl(\dfrac1{r}+\dfrac1{\theta_1}-\dfrac1m\Bigr)>0&\text{if }\ \dfrac1{r}+\dfrac1{\theta_1}\le1,\\[2mm]
s_1+n\Bigl(1-\dfrac1m\Bigr)>0&\text{if }\ \dfrac1{r}+\dfrac1{\theta_1}>1,
\end{cases}
\end{equation}
and
\begin{equation}\label{eq:comp-lowhigh}
s_1+n\Bigl(\frac1{r}+\frac1{\theta_2}-\frac1m\Bigr)<\nu .
\end{equation}
Then there is a constant $C_0>0$, depending on $n,p,s,\sigma,r,m,\theta_1,\theta_2$ and $\nu$, such that
\begin{equation}\label{eq:composition}
\|F(u)\|_{\B{s_1}m\sigma}\le C_0\|u\|_{\B sr\infty}^{p-1}\|u\|_{\B{s}{r}\sigma}
\end{equation}
for every $u\in\B{s}{r}\sigma$.
\end{lemma}

\begin{proof}
We write $M_k:=M_k(u)$ and
\[
c_k:=2^{ks}\|\Dj ku\|_{L^{r}},\qquad \|(c_k)_{k\in\Z}\|_{\ell^\sigma}=\|u\|_{\B{s}{r}\sigma}.
\]
By \eqref{eq:comp-structural}, Lemma \ref{lem:gain} (i) applies with $\theta=\theta_1$ and with $\theta=\theta_2$. Hence, for $i=1,2$ and all $k,l\in\Z$,
\begin{equation}\label{eq:Mk-coeff}
\|M_k\|_{L^{\theta_i}}\le C2^{k((p-1)\kappa-\frac n{\theta_i})}\|u\|_{\B sr\infty}^{p-1},\qquad
\|\Dj lM_k\|_{L^{\theta_2}}\le C2^{(k-l)\nu}2^{k((p-1)\kappa-\frac n{\theta_2})}\|u\|_{\B sr\infty}^{p-1}.
\end{equation}

From \eqref{eq:telescopic} and the identity $M_k=\Sj{k-3}M_k+\tDj kM_k+\sum_{l\ge k+3}\Dj lM_k$ in $\Sch'(\Rn)$ we obtain the paraproduct splitting
\begin{equation}\label{eq:decomposition}
\begin{aligned}
\Dj jF(u)&=\sum_{k\in\Z}\Dj j\bigl(M_k\Dj ku\bigr)\\
&=\sum_{k\in\Z}\Dj j\bigl(\Sj{k-3}(M_k)\Dj ku\bigr)+\sum_{k\in\Z}\Dj j\bigl(\tDj k(M_k)\Dj ku\bigr)+\sum_{k\in\Z}\sum_{l\ge k+3}\Dj j\bigl(\Dj lM_k\,\Dj ku\bigr)\\
&=:I_j+II_j+III_j .
\end{aligned}
\end{equation}
Below we define $\beta_1,\beta_2\in(0,\infty]$ by $\frac1{\beta_i}:=\frac1{r}+\frac1{\theta_i}$. By \eqref{eq:comp-structural}, $\frac1m\le\frac1{\beta_i}$. Moreover, \eqref{eq:s1} gives
\begin{equation}\label{eq:exponent-identity}
(p-1)\kappa-\frac n{\theta_i}-s=p\kappa-\frac n{\beta_i}=-\Bigl(s_1+\frac n{\beta_i}-\frac nm\Bigr).
\end{equation}

\emph{Estimate of $I_j$.} Since
\begin{equation}\label{eq:supp-I}
\supp\F\bigl[\Sj{k-3}(M_k)\Dj ku\bigr]\subset\{\xi\in\Rn:\ 2^{k-2}\le|\xi|\le2^{k+2}\},
\end{equation}
we have
\[
I_j=\sum_{|k-j|\le4}\Dj j\bigl(\Sj{k-3}(M_k)\Dj ku\bigr).
\]
Hence
\[
\begin{aligned}
\|I_j\|_{L^m}&\le C\sum_{|k-j|\le4}\bigl\|\Sj{k-3}(M_k)\Dj ku\bigr\|_{L^m}
\le C\sum_{|k-j|\le4}2^{kn(\frac1{\beta_1}-\frac1m)}\bigl\|\Sj{k-3}(M_k)\Dj ku\bigr\|_{L^{\beta_1}}\\
&\le C\sum_{|k-j|\le4}2^{kn(\frac1{\beta_1}-\frac1m)}\|M_k\|_{L^{\theta_1}}\|\Dj ku\|_{L^{r}},
\end{aligned}
\]
where we used \eqref{eq:supp-I}, Proposition \ref{prop:Bernstein} with $0<\beta_1\le m$, and the Hölder inequality, which also holds for exponents below one (see \cite[p.~10]{Gra08}). Using \eqref{eq:Mk-coeff} and \eqref{eq:exponent-identity}, we get
\[
\|I_j\|_{L^m}\le C\sum_{|k-j|\le4}2^{k(\frac n{\beta_1}-\frac nm+(p-1)\kappa-\frac n{\theta_1}-s)}c_k\|u\|_{\B sr\infty}^{p-1}=C\sum_{|k-j|\le4}2^{-ks_1}c_k\|u\|_{\B sr\infty}^{p-1},
\]
and therefore
\begin{equation}\label{eq:I-final}
\Bigl\|\bigl(2^{js_1}\|I_j\|_{L^m}\bigr)_{j\in\Z}\Bigr\|_{\ell^\sigma}\le C\|u\|_{\B sr\infty}^{p-1}\|u\|_{\B{s}{r}\sigma}
\end{equation}
for every $1\le\sigma\le\infty$.

\emph{Estimate of $II_j$.} Since
\begin{equation}\label{eq:supp-II}
\supp\F\bigl[\tDj k(M_k)\Dj ku\bigr]\subset\{\xi\in\Rn:\ |\xi|\le2^{k+4}\},
\end{equation}
we have
\[
II_j=\sum_{k\ge j-5}\Dj j\bigl(\tDj k(M_k)\Dj ku\bigr).
\]
Assume first that $\frac1{r}+\frac1{\theta_1}\le1$, that is, $\beta_1\ge1$. In this case $\Dj j$ and $\tDj k$ are bounded on $L^{\beta_1}(\Rn)$ and on $L^{\theta_1}(\Rn)$. Proposition \ref{prop:Bernstein} at frequency $2^j$ and \eqref{eq:Mk-coeff} give
\[
\begin{aligned}
\|II_j\|_{L^m}&\le C\sum_{k\ge j-5}2^{jn(\frac1{\beta_1}-\frac1m)}\bigl\|\tDj k(M_k)\Dj ku\bigr\|_{L^{\beta_1}}
\le C\sum_{k\ge j-5}2^{jn(\frac1{\beta_1}-\frac1m)}\|M_k\|_{L^{\theta_1}}\|\Dj ku\|_{L^{r}}\\
&\le C\sum_{k\ge j-5}2^{jn(\frac1{\beta_1}-\frac1m)}2^{k((p-1)\kappa-\frac n{\theta_1}-s)}c_k\|u\|_{\B sr\infty}^{p-1}.
\end{aligned}
\]
By \eqref{eq:exponent-identity}, this means
\[
2^{js_1}\|II_j\|_{L^m}\le C\sum_{k\ge j-5}2^{(j-k)(s_1+\frac n{\beta_1}-\frac nm)}c_k\|u\|_{\B sr\infty}^{p-1}.
\]
The first condition in \eqref{eq:comp-highhigh} and the Young inequality for convolutions then give
\begin{equation}\label{eq:II-final}
\Bigl\|\bigl(2^{js_1}\|II_j\|_{L^m}\bigr)_{j\in\Z}\Bigr\|_{\ell^\sigma}\le C\|u\|_{\B sr\infty}^{p-1}\|u\|_{\B{s}{r}\sigma}
\end{equation}
for every $1\le\sigma\le\infty$. Assume now that $\frac1{r}+\frac1{\theta_1}>1$, that is, $\beta_1<1$. Then
\[
\begin{aligned}
\|II_j\|_{L^m}&\le C\sum_{k\ge j-5}2^{jn(1-\frac1m)}\bigl\|\tDj k(M_k)\Dj ku\bigr\|_{L^1}
\le C2^{jn(1-\frac1m)}\sum_{k\ge j-5}2^{kn(\frac1{\beta_1}-1)}\bigl\|\tDj k(M_k)\Dj ku\bigr\|_{L^{\beta_1}}\\
&\le C2^{jn(1-\frac1m)}\sum_{k\ge j-5}2^{kn(\frac1{\beta_1}-1)}\|M_k\|_{L^{\theta_1}}\|\Dj ku\|_{L^{r}},
\end{aligned}
\]
where we applied Proposition \ref{prop:Bernstein} at frequency $2^j$ with $1\le m$, and at frequency $2^{k+4}$ with $\beta_1<1$, which is possible by \eqref{eq:supp-II}. By \eqref{eq:exponent-identity}, $\frac n{\beta_1}-n+(p-1)\kappa-\frac n{\theta_1}-s=-\bigl(s_1+n-\frac nm\bigr)$, and \eqref{eq:Mk-coeff} gives
\[
2^{js_1}\|II_j\|_{L^m}\le C\sum_{k\ge j-5}2^{(j-k)(s_1+n-\frac nm)}c_k\|u\|_{\B sr\infty}^{p-1}.
\]
The second condition in \eqref{eq:comp-highhigh} now gives \eqref{eq:II-final}.

\emph{Estimate of $III_j$.} Since
\begin{equation}\label{eq:supp-III}
\supp\F\bigl[\Dj l(M_k)\Dj ku\bigr]\subset\{\xi\in\Rn:\ 2^{l-2}\le|\xi|\le2^{l+2}\}\qquad\text{for }l\ge k+3,
\end{equation}
we have
\[
III_j=\sum_{k\in\Z}\ \sum_{\substack{l\ge k+3\\|l-j|\le4}}\Dj j\bigl(\Dj l(M_k)\Dj ku\bigr).
\]
In this sum $k\le l-3\le j+1$. Hence
\[
\begin{aligned}
\|III_j\|_{L^m}&\le C\sum_{|l-j|\le4}\ \sum_{k\le l-3}2^{ln(\frac1{\beta_2}-\frac1m)}\|\Dj lM_k\|_{L^{\theta_2}}\|\Dj ku\|_{L^{r}}\\
&\le C\sum_{|l-j|\le4}\ \sum_{k\le l-3}2^{ln(\frac1{\beta_2}-\frac1m)}2^{(k-l)\nu}2^{k((p-1)\kappa-\frac n{\theta_2}-s)}c_k\|u\|^{p-1}\\
&\le C\sum_{k\le j+1}2^{j(\frac n{\beta_2}-\frac nm-\nu)}2^{k(\nu+(p-1)\kappa-\frac n{\theta_2}-s)}c_k\|u\|_{\B sr\infty}^{p-1},
\end{aligned}
\]
where we used \eqref{eq:supp-III}, Proposition \ref{prop:Bernstein} with $0<\beta_2\le m$, the Hölder inequality and the second bound in \eqref{eq:Mk-coeff}. The gain $2^{(k-l)\nu}$ is the key point here. By \eqref{eq:exponent-identity}, $\nu+(p-1)\kappa-\frac n{\theta_2}-s=-\bigl(s_1+\frac n{\beta_2}-\frac nm-\nu\bigr)$, so that
\[
2^{js_1}\|III_j\|_{L^m}\le C\sum_{k\le j+1}2^{(j-k)(s_1+\frac n{\beta_2}-\frac nm-\nu)}c_k\|u\|_{\B sr\infty}^{p-1}.
\]
The convolution inequality then gives
\begin{equation}\label{eq:III-final}
\Bigl\|\bigl(2^{js_1}\|III_j\|_{L^m}\bigr)_{j\in\Z}\Bigr\|_{\ell^\sigma}\le C\|u\|_{\B sr\infty}^{p-1}\|u\|_{\B{s}{r}\sigma},
\end{equation}
provided that $s_1+\frac n{\beta_2}-\frac nm-\nu<0$, which is \eqref{eq:comp-lowhigh}.

The lemma follows from \eqref{eq:I-final}, \eqref{eq:II-final} and \eqref{eq:III-final}.
\end{proof}

\begin{remark}\label{rem:Chikami-recovered}
For $\nu=0$, Lemma \ref{lem:composition} reduces to the non-endpoint part of \cite[Lemma 3.2]{Chi19}, with $p-1$, $r$ and $m$ in place of $b$, $q$ and $p$ in the notation used there. Indeed, in that case \eqref{eq:comp-structural} means $r\le(p-1)\theta_i$, $\frac1m\le\frac1r+\frac1{\theta_i}$ and $s<n(\frac1r-\frac1{(p-1)\theta_i})$. Condition \eqref{eq:comp-highhigh} reads $s>\frac{(p-1)n}p(\frac1r-\frac1{(p-1)\theta_1})$ if $\frac1r+\frac1{\theta_1}\le1$, and $s>n(\frac1r-\frac1p)$ if $\frac1r+\frac1{\theta_1}>1$. Condition \eqref{eq:comp-lowhigh} reads $s<\frac{(p-1)n}p(\frac1r-\frac1{(p-1)\theta_2})$. These are exactly the conditions of the parabolic estimate. Observe that, for $\nu=0$, \eqref{eq:comp-lowhigh} forces $s_1<0$, so the nonlinearity is measured in a Besov space of negative regularity. The gain $\nu>0$ is the new point needed for the wave equation, where $F(u)$ has to be measured in $\B{\sigma_r}{r'}\infty$ with $\sigma_r>0$. 
\end{remark}

\subsection{Difference estimate}

To close the fixed point argument we also need a difference estimate. Once again, Proposition \ref{prop:Bernstein} is the main tool.

\begin{lemma}\label{lem:difference}
Let $n\ge1$, $F(u):=|u|^{p-1}u$ with $\R\ni p\ge2$, $1\le\sigma\le\infty$ and $0\le\nu\le p-2$. Assume that $s,r,m,\theta_1,\theta_2$ satisfy \eqref{eq:comp-structural}, \eqref{eq:comp-highhigh} and \eqref{eq:comp-lowhigh}, with $s_1$ given by \eqref{eq:s1}. Then there is a constant $C_1>0$, depending on $n,p,s,\sigma,r,m,\theta_1,\theta_2$ and $\nu$, such that
\begin{equation}\label{eq:difference}
\begin{aligned}
\|F(u)-F(v)\|_{\B{s_1}m\sigma}\le C_1\Bigl\{&\bigl(\|u\|_{\B sr\infty}+\|v\|_{\B sr\infty}\bigr)^{p-2}\|u-v\|_{\B sr\infty}\|u\|_{\B{s}{r}\sigma}\\
&+\|v\|_{\B sr\infty}^{p-1}\|u-v\|_{\B{s}{r}\sigma}\Bigr\}
\end{aligned}
\end{equation}
for all $u,v\in\B{s}{r}\sigma$. In particular, for $\sigma=\infty$,
\begin{equation}\label{eq:difference-infty}
\|F(u)-F(v)\|_{\B{s_1}m\infty}\le 2C_1\bigl(\|u\|_{\B sr\infty}+\|v\|_{\B sr\infty}\bigr)^{p-1}\|u-v\|_{\B sr\infty}.
\end{equation}
\end{lemma}

\begin{proof}
Representing the power function by the telescopic series \eqref{eq:telescopic} for both $u$ and $v$, we have
\[
F(u)-F(v)=\sum_{k\in\Z}\bigl(M_k(u)-M_k(v)\bigr)\Dj ku+\sum_{k\in\Z}M_k(v)\Dj k(u-v)=:J_1+J_2 .
\]
Put $w=u-v$. The term
\[
J_2=\sum_{k\in\Z}M_k(v)\Dj kw
\]
is treated exactly like $F(u)$ in the proof of Lemma \ref{lem:composition}, with $M_k(v)$ in place of $M_k(u)$ and $\Dj kw$ in place of $\Dj ku$. This is possible because $\nu\le p-2\le p-1$. Under our assumptions on $s,\sigma,r,m,\theta_1,\theta_2$ and $\nu$ we thus have
\begin{equation}\label{eq:J2}
\Bigl\|\bigl(2^{js_1}\|\Dj jJ_2\|_{L^m}\bigr)_{j\in\Z}\Bigr\|_{\ell^\sigma}\le C\|v\|_{\B sr\infty}^{p-1}\|w\|_{\B{s}{r}\sigma}.
\end{equation}

We bound $J_1$. We set $\widetilde M_k:=M_k(u)-M_k(v)$ and $\|\cdot\|:=\|\cdot\|_{\B sr\infty}$. The mean value formula \eqref{eq:F3} gives the representation \eqref{eq:tildeM-repr}, namely
\[
\widetilde M_k=\int_0^1\bigl(\Sj{k-1}w+t\Dj kw\bigr)\int_0^1F''\bigl(\lambda z^t_k(u)+(1-\lambda)z^t_k(v)\bigr)d\lambda\,dt .
\]
By \eqref{eq:comp-structural}, Lemma \ref{lem:gain} (ii) applies with $\theta=\theta_1$ and with $\theta=\theta_2$. Hence, for $i=1,2$ and all $k,l\in\Z$,
\begin{equation}\label{eq:tildeM-coeff}
\begin{aligned}
\|\widetilde M_k\|_{L^{\theta_i}}&\le C2^{k((p-1)\kappa-\frac n{\theta_i})}\bigl(\|u\|+\|v\|\bigr)^{p-2}\|w\|,\\
\|\Dj l\widetilde M_k\|_{L^{\theta_2}}&\le C2^{(k-l)\nu}2^{k((p-1)\kappa-\frac n{\theta_2})}\bigl(\|u\|+\|v\|\bigr)^{p-2}\|w\|.
\end{aligned}
\end{equation}

 The Bony decomposition gives
\[
\begin{aligned}
\Dj jJ_1&=\sum_{|k-j|\le4}\Dj j\bigl(\Sj{k-3}(\widetilde M_k)\Dj ku\bigr)+\sum_{k\ge j-5}\Dj j\bigl(\tDj k(\widetilde M_k)\Dj ku\bigr)
+\sum_{k\in\Z}\Dj j\Bigl(\sum_{\substack{l\ge k+3\\|l-j|\le4}}\Dj l\widetilde M_k\,\Dj ku\Bigr)\\
&=:I_j+II_j+III_j .
\end{aligned}
\]
The bounds \eqref{eq:tildeM-coeff} coincide with \eqref{eq:Mk-coeff}, except that $\|u\|^{p-1}$ is replaced by $(\|u\|+\|v\|)^{p-2}\|w\|$. Hence $I_j$, $II_j$ and $III_j$ are estimated as in the proof of Lemma \ref{lem:composition}, and we omit some details. Let $c_k:=2^{ks}\|\Dj ku\|_{L^{r}}$ and $\frac1{\beta_i}=\frac1{r}+\frac1{\theta_i}$.

\emph{Estimate of $I_j$.} As for $I_j$ in the proof of Lemma \ref{lem:composition}, Proposition \ref{prop:Bernstein}, the Hölder inequality, \eqref{eq:exponent-identity} and \eqref{eq:tildeM-coeff} give
\[
2^{js_1}\|I_j\|_{L^m}\le C\sum_{|k-j|\le4}2^{(j-k)s_1}c_k\,\bigl(\|u\|+\|v\|\bigr)^{p-2}\|w\|,
\]
so that
\[
\Bigl\|\bigl(2^{js_1}\|I_j\|_{L^m}\bigr)_{j\in\Z}\Bigr\|_{\ell^\sigma}\le C\bigl(\|u\|_{\B sr\infty}+\|v\|_{\B sr\infty}\bigr)^{p-2}\|w\|_{\B sr\infty}\|u\|_{\B{s}{r}\sigma}
\]
for $1\le\sigma\le\infty$.

\emph{Estimate of $II_j$.} We distinguish the cases $\frac1{r}+\frac1{\theta_1}\le1$ and $\frac1{r}+\frac1{\theta_1}>1$, exactly as in the proof of Lemma \ref{lem:composition}. Using \eqref{eq:tildeM-coeff}, we obtain
\[
2^{js_1}\|II_j\|_{L^m}\le
\begin{cases}
C\displaystyle\sum_{k\ge j-5}2^{(j-k)(s_1+\frac n{\beta_1}-\frac nm)}c_k\,\bigl(\|u\|+\|v\|\bigr)^{p-2}\|w\|&\text{if }\beta_1\ge1,\\[3mm]
C\displaystyle\sum_{k\ge j-5}2^{(j-k)(s_1+n-\frac nm)}c_k\,\bigl(\|u\|+\|v\|\bigr)^{p-2}\|w\|&\text{if }\beta_1<1.
\end{cases}
\]
Then \eqref{eq:comp-highhigh} and the Young inequality give
\[
\Bigl\|\bigl(2^{js_1}\|II_j\|_{L^m}\bigr)_{j\in\Z}\Bigr\|_{\ell^\sigma}\le C\bigl(\|u\|_{\B sr\infty}+\|v\|_{\B sr\infty}\bigr)^{p-2}\|w\|_{\B sr\infty}\|u\|_{\B{s}{r}\sigma}
\]
for all $1\le\sigma\le\infty$.

\emph{Estimate of $III_j$.} By \eqref{eq:supp-III}, Proposition \ref{prop:Bernstein}, the Hölder inequality and the second bound in \eqref{eq:tildeM-coeff}, we have
\[
2^{js_1}\|III_j\|_{L^m}\le C\sum_{k\le j+1}2^{(j-k)(s_1+\frac n{\beta_2}-\frac nm-\nu)}c_k\,\bigl(\|u\|+\|v\|\bigr)^{p-2}\|w\|,
\]
and \eqref{eq:comp-lowhigh} gives
\[
\Bigl\|\bigl(2^{js_1}\|III_j\|_{L^m}\bigr)_{j\in\Z}\Bigr\|_{\ell^\sigma}\le C\bigl(\|u\|_{\B sr\infty}+\|v\|_{\B sr\infty}\bigr)^{p-2}\|w\|_{\B sr\infty}\|u\|_{\B{s}{r}\sigma}.
\]
Combining these bounds with \eqref{eq:J2}, we obtain \eqref{eq:difference}. Finally, \eqref{eq:difference-infty} follows from \eqref{eq:difference} with $\sigma=\infty$, because $\|u\|^{p-2}_{\B sr\infty}\|u\|_{\B sr\infty}+\|v\|^{p-1}_{\B sr\infty}\le2(\|u\|_{\B sr\infty}+\|v\|_{\B sr\infty})^{p-1}$.
\end{proof}
\begin{remark}\label{rem:diff}
Unlike Lemma \ref{lem:composition}, Lemma \ref{lem:difference} requires $p\ge2$, because the proof uses \eqref{eq:F3} through Lemma \ref{lem:gain} (ii). The gain is now limited by the smoothness of $F''$, which requires $\nu\le p-2$. For $\nu=0$, the right-hand side of \eqref{eq:difference} is bounded by $C(\|u\|_{\B sr\infty}^{p-1}+\|v\|_{\B sr\infty}^{p-1})\|u-v\|_{\B sr\sigma}$, which is the non-endpoint part of \cite[Lemma 3.7]{Chi19}.
\end{remark}

\subsection{Application to the Hardy--Hénon nonlinearity}

We now put together the composition and difference estimates and the estimate for the singular weight. The next proposition applies the previous estimates to the weighted nonlinearity $h_\gamma F(u)$.

\begin{proposition}\label{prop:weighted-nonlinear}
Let $n\ge2$, $0<\gamma<\min\{2,n\}$, $p>2$ and $2<r<\infty$ satisfy \eqref{eq:H2} and \eqref{eq:H3}, and let $s_r$, $\sigma_r$ be given by \eqref{eq:indices}. Then $0<s_r<\frac nr$ and $r<p+1$. Moreover, there is a constant $C>0$ such that
\begin{align}
\|F(u)\|_{\B{\sigma_r}{r'}\infty}&\le C\|u\|_{\B{s_r}r\infty}^p,\label{eq:unweighted-comp}\\
\|F(u)-F(v)\|_{\B{\sigma_r}{r'}\infty}&\le C\bigl(\|u\|_{\B{s_r}r\infty}+\|v\|_{\B{s_r}r\infty}\bigr)^{p-1}\|u-v\|_{\B{s_r}r\infty},\label{eq:unweighted-diff}\\
\|h_\gamma F(u)\|_{\B{s_r-1+2\beta_r}{r'}\infty}&\le C\|u\|_{\B{s_r}r\infty}^p,\label{eq:weighted-comp}\\
\|h_\gamma(F(u)-F(v))\|_{\B{s_r-1+2\beta_r}{r'}\infty}&\le C\bigl(\|u\|_{\B{s_r}r\infty}+\|v\|_{\B{s_r}r\infty}\bigr)^{p-1}\|u-v\|_{\B{s_r}r\infty},\label{eq:weighted-diff}
\end{align}
for all $u,v\in\B{s_r}r\infty(\Rn)$. 
\end{proposition}

\begin{proof}
We use Lemmas \ref{lem:composition} and \ref{lem:difference} with
\[
s:=s_r,\qquad m:=r',\qquad \theta_1=\theta_2:=\frac r{r-2},\qquad \sigma:=\infty,\qquad \nu:=p-2,
\]
so that $\frac1r+\frac1{\theta_i}=1-\frac1r=\frac1m$, $\beta_1=\beta_2=m$ and $\nu=p-2\le p-1$. By \eqref{eq:homogeneity}, the index $s_1$ of \eqref{eq:s1} is $\frac n{r'}-p(\frac nr-s_r)=\sigma_r$. Conditions \eqref{eq:comp-highhigh} and \eqref{eq:comp-lowhigh} therefore read $0<\sigma_r<\nu$, which is part of \eqref{eq:H3}. We now verify \eqref{eq:comp-structural}. With $\kappa=\frac nr-s_r$, a direct computation from \eqref{eq:beta-eta} and \eqref{eq:indices} shows that
\[
(p-1)\kappa-\frac n{\theta_i}=1-\gamma+\eta_r-n\Bigl(1-\frac2r\Bigr)=1-\gamma-2\beta_r,
\]
because $n(1-\frac2r)-\eta_r=(n+1)(\frac12-\frac1r)=2\beta_r$. Hence $\frac n{\theta_i}<(p-1)\kappa$ is exactly \eqref{eq:H2}, and in particular $\kappa>0$. The last condition in \eqref{eq:comp-structural} holds with equality. We still have to check $r\le(p-1)\theta_i$, which means $r\le p+1$. By \eqref{eq:s1}, $\sigma_r>0$ is equivalent to $s_r>\frac{n(p+1-r)}{pr}$. On the other hand, $\frac n{\theta_i}<(p-1)\kappa$ is equivalent to $s_r<\frac{n(p+1-r)}{(p-1)r}$. These two inequalities together give $(p+1-r)(\frac1{p-1}-\frac1p)>0$, that is, $r<p+1$. Consequently $s_r>\frac{n(p+1-r)}{pr}>0$.  If $r>p$, then $b=r$ and $\sum_{j>0}\|\Dj ju\|_{L^r}\le\sum_{j>0}2^{-js_r}\|u\|_{\B{s_r}r\infty}<\infty$, since $s_r>0$. If $r\le p$, then $b=p$ and Proposition \ref{prop:Bernstein} gives
\[
\sum_{j>0}\|\Dj ju\|_{L^p}\le C\sum_{j>0}2^{-j(s_r-\frac nr+\frac np)}\|u\|_{\B{s_r}r\infty}<\infty,
\]
because $s_r-\frac nr+\frac np>\frac{n(p+1-r)}{pr}-\frac nr+\frac np=\frac n{pr}>0$. Lemmas \ref{lem:composition} and \ref{lem:difference} therefore give \eqref{eq:unweighted-comp} and \eqref{eq:unweighted-diff}.
It remains to remove the weight. We apply Lemma \ref{lem:potential} with $m=r'$ and $\varsigma=\sigma_r$, which requires $\gamma-\frac nr<\sigma_r<\frac n{r'}$. The first inequality is part of \eqref{eq:H3}. The second one holds because $\sigma_r=\frac n{r'}-p\kappa$ by \eqref{eq:homogeneity}, and $\kappa>0$. Since $\sigma_r-\gamma=s_r-1+2\beta_r$, this proves \eqref{eq:weighted-comp}. Applying Lemma \ref{lem:potential} to $F(u)-F(v)$ gives \eqref{eq:weighted-diff}.
\end{proof}

\begin{remark}\label{rem:weighted-nonlinear}
The choice $\theta_1=\theta_2=\frac r{r-2}$ in the proof above comes from the fact that the output exponent $r'$ is smaller than the input exponent $r$. The coefficients $M_k(u)$ must then be estimated in $L^{\theta}$ with $\frac1\theta\ge\frac1{r'}-\frac1r$, and the smallest admissible value of $\frac1\theta$ gives the weakest form of \eqref{eq:comp-lowhigh}. With this choice, \eqref{eq:comp-structural} becomes \eqref{eq:H2}, and \eqref{eq:comp-highhigh}, \eqref{eq:comp-lowhigh} become the window $0<\sigma_r<\nu$ in \eqref{eq:H3}. The remaining condition $\sigma_r>\gamma-\frac nr$ in \eqref{eq:H3} comes from the product law of Abidi.
\end{remark}

\section{Global well-posedness in the time-weighted framework}\label{sec:wellposedness}

In this section the hypotheses of Theorem \ref{thm:global} are in force, and $s:=s_r$, $\zeta:=\zeta_r$. By Proposition \ref{prop:weighted-nonlinear}, $0<s<\frac nr$. We write $\kappa:=\frac nr-s$. Recall that \eqref{eq:integral} reads
\[
u(t)=u_L(t)+B(u)(t),\qquad B(u)(t)=a\int_0^tK(t-\tau)\bigl(h_\gamma F(u(\tau))\bigr)d\tau .
\]

\subsection{Estimates for the Duhamel term}

The next lemma controls the nonlinear term. It follows from the weighted dispersive estimate of Lemma \ref{lem:group-potential} and from Proposition \ref{prop:weighted-nonlinear}. The time weights then lead to a beta integral, which is finite because of \eqref{eq:H1} and \eqref{eq:H2}.

\begin{lemma}\label{lem:key}
There is a constant $C_0>0$, depending on $n,\gamma,p,r$, such that
\begin{align}
\|B(u)\|_{\E}&\le C_0|a|\,\|u\|_{\E}^p,\label{eq:key-1}\\
\|B(u)-B(v)\|_{\E}&\le C_0|a|\,\bigl(\|u\|_{\E}+\|v\|_{\E}\bigr)^{p-1}\|u-v\|_{\E},\label{eq:key-2}
\end{align}
for all $u,v\in\E$.
\end{lemma}

\begin{proof}
Let $t>0$. By the proof of Proposition \ref{prop:weighted-nonlinear}, $\gamma-\frac nr<\sigma_r<\frac n{r'}$, and $\sigma_r-\gamma+1-2\beta_r=s$. Using Lemma \ref{lem:group-potential} with $\theta=1$, $\sigma=\infty$ and $\varsigma=\sigma_r$, and \eqref{eq:unweighted-comp} therefore give
\[
\begin{aligned}
\|B(u)(t)\|_{\B sr\infty}&\le|a|\int_0^t\bigl\|K(t-\tau)\bigl(h_\gamma F(u(\tau))\bigr)\bigr\|_{\B sr\infty}d\tau\\
&\le C|a|\int_0^t(t-\tau)^{-\eta_r}\bigl\|F(u(\tau))\bigr\|_{\B{\sigma_r}{r'}\infty}d\tau\\
&\le C|a|\int_0^t(t-\tau)^{-\eta_r}\|u(\tau)\|_{\B sr\infty}^pd\tau .
\end{aligned}
\]
By definition of $\E$, the right-hand side is at most
\[
C|a|\,\|u\|_{\E}^p\int_0^t(t-\tau)^{-\eta_r}\tau^{-p\zeta}d\tau
=C|a|\,\|u\|_{\E}^p\,t^{1-\eta_r-p\zeta}\int_0^1(1-z)^{-\eta_r}z^{-p\zeta}dz
=C_0|a|\,\|u\|_{\E}^p\,t^{-\zeta}.
\]
Here we used that $\eta_r<1$ and $p\zeta<1$ (Remark \ref{rem:hypotheses} (a)), so that the beta integral is finite, and that $1-\eta_r-p\zeta=-\zeta$ by \eqref{eq:zeta}. Multiplying by $t^\zeta$ and taking the supremum over $t>0$ gives \eqref{eq:key-1}. In the same way, Lemma \ref{lem:group-potential} and \eqref{eq:unweighted-diff} give
\[
\begin{aligned}
\|B(u)(t)-B(v)(t)\|_{\B sr\infty}
&\le C|a|\int_0^t(t-\tau)^{-\eta_r}\bigl\|F(u(\tau))-F(v(\tau))\bigr\|_{\B{\sigma_r}{r'}\infty}d\tau\\
&\le C|a|\int_0^t(t-\tau)^{-\eta_r}\bigl(\|u(\tau)\|_{\B sr\infty}+\|v(\tau)\|_{\B sr\infty}\bigr)^{p-1}\|u(\tau)-v(\tau)\|_{\B sr\infty}d\tau\\
&\le C|a|\bigl(\|u\|_{\E}+\|v\|_{\E}\bigr)^{p-1}\|u-v\|_{\E}\int_0^t(t-\tau)^{-\eta_r}\tau^{-p\zeta}d\tau\\
&=C_0|a|\bigl(\|u\|_{\E}+\|v\|_{\E}\bigr)^{p-1}\|u-v\|_{\E}\,t^{-\zeta},
\end{aligned}
\]
which proves \eqref{eq:key-2}.
\end{proof}

\subsection{Global well-posedness and self-similarity}

With Lemma \ref{lem:key} at hand, we construct global solutions by a contraction argument in a closed ball of $\E$. The next theorem gives existence, uniqueness and Lipschitz dependence on the data when the data are small in $\Ic$. It also shows that homogeneous data of the degrees given by the scaling \eqref{eq:scaling} produce self-similar solutions.

\begin{theorem}[Global well-posedness and self-similarity]\label{thm:global}
Let $n\ge2$, $a\in\R$, $0<\gamma<\min\{2,n\}$, $p>2$ and $2<r<\infty$, and assume that
\begin{align}
&\frac1p<\eta_r,\label{eq:H1}\\
&2\beta_r+\gamma<1,\label{eq:H2}\\
&\max\Bigl\{0,\gamma-\frac nr\Bigr\}<\sigma_r<p-2.\label{eq:H3}
\end{align}
Consider $s=s_r$ and $\zeta=\zeta_r$ defined by \eqref{eq:indices}. Then there exists $\varepsilon>0$ such that the following statements hold.
\begin{enumerate}[(i)]
\item (Well-posedness) If $(u_0,u_1)\in\Ic$ satisfies
\begin{equation}\label{eq:small}
\|(u_0,u_1)\|_{\Ic}\le\varepsilon,
\end{equation}
then \eqref{eq:integral} has a unique solution $u\in\E$ such that $\|u\|_{\E}\le2\varepsilon$. Moreover, if $v$ is the solution associated with $(v_0,v_1)\in\Ic$, where $\|(v_0,v_1)\|_{\Ic}\le\varepsilon$, then
\[
\|u-v\|_{\E}\le2\|(u_0,u_1)-(v_0,v_1)\|_{\Ic}.
\]
\item (Self-similarity) If, in addition, $u_0$ and $u_1$ are homogeneous of degrees $-\frac{2-\gamma}{p-1}$ and $-\frac{2-\gamma}{p-1}-1$, respectively, then the solution $u$  of \eqref{eq:integral} is self-similar.
\end{enumerate}
\end{theorem}

\begin{proof}
If $a=0$ then $u=u_L$, so we may assume $a\ne0$. Let $C_0$ be the constant of Lemma \ref{lem:key}, and choose $\varepsilon>0$ such that
\begin{equation}\label{eq:eps}
2\cdot4^{p-1}C_0|a|\,\varepsilon^{p-1}\le1 .
\end{equation}
Let $(u_0,u_1)\in\Ic$ satisfy \eqref{eq:small}. We consider the closed ball $\overline B(0,2\varepsilon)=\{u\in\E:\ \|u\|_{\E}\le2\varepsilon\}$ and the map
\[
\Phi(u)(t):=u_L(t)+B(u)(t),\qquad t>0 .
\]
We show that $\Phi$ is a contraction of $\overline B(0,2\varepsilon)$ into itself. For $u\in\overline B(0,2\varepsilon)$, \eqref{eq:key-1}, \eqref{eq:small} and \eqref{eq:eps} give
\[
\|\Phi(u)\|_{\E}\le\|u_L\|_{\E}+\|B(u)\|_{\E}\le\varepsilon+C_0|a|(2\varepsilon)^p\le\varepsilon+\varepsilon=2\varepsilon,
\]
because $C_0|a|2^p\varepsilon^{p-1}\le2\cdot4^{p-1}C_0|a|\varepsilon^{p-1}\le1$. So $\Phi$ maps the ball into itself. For $u,v\in\overline B(0,2\varepsilon)$, \eqref{eq:key-2} and \eqref{eq:eps} give
\[
\|\Phi(u)-\Phi(v)\|_{\E}=\|B(u)-B(v)\|_{\E}\le C_0|a|(4\varepsilon)^{p-1}\|u-v\|_{\E}\le\frac12\|u-v\|_{\E}.
\]
The Banach fixed point theorem provides a unique fixed point $u$ of $\Phi$ in $\overline B(0,2\varepsilon)$. This is the unique solution of \eqref{eq:integral} in $\E$ with $\|u\|_{\E}\le2\varepsilon$. If $v$ is the solution with data $(v_0,v_1)$ and linear part $v_L$, then
\[
\|u-v\|_{\E}\le\|u_L-v_L\|_{\E}+\|B(u)-B(v)\|_{\E}\le\|(u_0,u_1)-(v_0,v_1)\|_{\Ic}+\frac12\|u-v\|_{\E},
\]
since $(u_0,u_1)\mapsto u_L$ is linear. This is the Lipschitz dependence, and (i) is proved.

We turn to (ii). Assume that $u_0$ and $u_1$ are homogeneous of degrees $-\frac{2-\gamma}{p-1}$ and $-\frac{2-\gamma}{p-1}-1$, and let $u$ be the solution from (i). Consider the Picard iterates
\[
u^1:=u_L,\qquad u^{m+1}:=u_L+B(u^m),\qquad m\in\mathbb N.
\]
The contraction argument shows that $u^m\in\overline B(0,2\varepsilon)$ for every $m$ and that $u^m\to u$ in $\E$. Recall the notation $u_\lambda$ from \eqref{eq:scaling}, and put $\alpha_0:=\frac{2-\gamma}{p-1}$. The homogeneity of the data means $\lambda^{\alpha_0}u_0(\lambda\cdot)=u_0$ and $\lambda^{\alpha_0+1}u_1(\lambda\cdot)=u_1$. Moreover, $K(t)[f(\lambda\cdot)]=\lambda^{-1}(K(\lambda t)f)(\lambda\cdot)$ and $\dK(t)[f(\lambda\cdot)]=(\dK(\lambda t)f)(\lambda\cdot)$, and therefore $(u_L)_\lambda=u_L$. Since $h_\gamma(x)\,g(\lambda x)=\lambda^\gamma(h_\gamma g)(\lambda x)$ and $\alpha_0(p-1)=2-\gamma$, the change of variables $\tau'=\lambda\tau$ gives
\[
B(w_\lambda)(t)=a\lambda^{\alpha_0p+\gamma-1}\int_0^t\bigl(K(\lambda(t-\tau))h_\gamma F(w(\lambda\tau))\bigr)(\lambda\cdot)\,d\tau=\lambda^{\alpha_0}B(w)(\lambda t)(\lambda\cdot)=\bigl(B(w)\bigr)_\lambda(t)
\]
for every $w\in\E$. An induction on $m$ then shows that $u^m_\lambda=u^m$ for all $m\in\mathbb N$ and $\lambda>0$. On the other hand, $\|f(\lambda\cdot)\|_{\B sr\infty}\le C\lambda^{s-\frac nr}\|f\|_{\B sr\infty}$ with $C$ independent of $\lambda$ (see \cite[Chapter 2]{BCD11}), and $\zeta=\alpha_0+s-\frac nr$ by \eqref{eq:zeta}. Hence $\|w_\lambda\|_{\E}\le C\|w\|_{\E}$ for all $w\in\E$ and $\lambda>0$. It follows that
\[
\|u-u_\lambda\|_{\E}\le\|u-u^m\|_{\E}+\|u^m_\lambda-u_\lambda\|_{\E}\le(1+C)\|u-u^m\|_{\E}\to0\quad\text{as }m\to\infty .
\]
Thus $u=u_\lambda$, and $u$ is self-similar.
\end{proof}

\begin{remark}\label{rem:hypotheses}
We comment on the hypotheses of Theorem \ref{thm:global}.
\begin{enumerate}[(a)]
\item Since $\eta_r<2\beta_r$ by \eqref{eq:beta-eta}, condition \eqref{eq:H2} gives $\eta_r<1$, so that $\zeta_r>0$. By \eqref{eq:zeta}, condition \eqref{eq:H1} is then equivalent to $p\zeta_r<1$. Hence \eqref{eq:H1} and \eqref{eq:H2} contain the two time-integrability conditions $\eta_r<1$ and $p\zeta_r<1$ of the case $\gamma=0$, plus the extra requirement $2\beta_r<1-\gamma$. For $\gamma=0$, \eqref{eq:H1} and \eqref{eq:H2} become the conditions $\frac1{\alpha+1}<\gamma_p$ and $p<\frac{2(n+1)}{n-1}$ of \cite[Theorem 3.10]{CF25}, written in the notation used there, where $\alpha=p-1$, $\gamma_p=\eta_r$ and $p$ stands for our $r$. In that case $s_r$ and $\zeta_r$ are the indices $s$ and $\zeta$ of that theorem. Condition \eqref{eq:H3} is what we need in order to apply the composition and difference estimates of Section \ref{sec:composition} to the weighted nonlinearity (Proposition \ref{prop:weighted-nonlinear}).
\item By \eqref{eq:H2}, $\gamma<1$ and $r<\frac{2(n+1)}{n-1}$. Together with \eqref{eq:H1}, it also gives $\frac1p<\eta_r<\frac{(n-1)(1-\gamma)}{n+1}$, hence $p>\frac{n+1}{(n-1)(1-\gamma)}$.
\item The hypotheses can be satisfied in every dimension. To see this, put $x:=\frac12-\frac1r$. A direct computation from \eqref{eq:indices} gives
\begin{equation}\label{eq:sigma-formula}
(p-1)\sigma_r=\frac{(p-1)(n-2)}2+p\gamma-1+(p-n)x .
\end{equation}
Let $n\ge3$, $0<\gamma<\frac12$, $p:=n+2$ and $r:=\frac{2(n+1)}n$. Then $x=\frac1{2(n+1)}$, $\eta_r=\frac{n-1}{2(n+1)}$ and $2\beta_r=\frac12$. Condition \eqref{eq:H2} holds because $\gamma<\frac12$, and \eqref{eq:H1} holds because $(n+2)(n-1)>2(n+1)$ for $n\ge3$. By \eqref{eq:sigma-formula},
\[
(n+1)\sigma_r=\frac{(n+1)(n-2)}2+(n+2)\gamma-1+\frac1{n+1}.
\]
The right-hand side is positive. Since $\gamma<\frac12$, it is smaller than $\frac{n^2}2$, which is smaller than $n(n+1)=(p-1)(p-2)$. Moreover, $\frac nr=\frac{n^2}{2(n+1)}>1>\gamma$. Hence \eqref{eq:H3} holds as well.
\item The regularity $\sigma_r$ grows with the dimension. Indeed, \eqref{eq:H2} gives $0<x<\frac1{n+1}$, hence $(p-n)x>-1$, and \eqref{eq:sigma-formula} together with $\gamma>0$ gives
\[
\sigma_r>\frac{n-2}2-\frac2{p-1}.
\]
Since $\sigma_r<p-2$ by \eqref{eq:H3}, the power must satisfy $p-2+\frac2{p-1}>\frac{n-2}2$. For $n=6$ this forces $p>3$ and $\sigma_r>1$, because $2-\frac2{p-1}\ge p-2$ for $2<p\le3$. Hence, for $n\ge6$, the gain $\nu>1$ in Lemmas \ref{lem:composition} and \ref{lem:difference} is indispensable, and for large $n$ gains of order $\frac n2$ are needed.
\end{enumerate}
\end{remark}

\subsection{Besov data and self-similar solutions}

The class $\Ic$ is defined through the free evolution $u_L$, so it is natural to ask which data belong to it. The class $\Ic$ contains critical Besov spaces and homogeneous data, as the next result shows. The relevant Besov spaces are built on $L^{r'}$, and their regularity is the one that makes their norm invariant under \eqref{eq:scaling}. The proof combines Lemma \ref{lem:group}, used with an intermediate value of $\theta$, and Proposition \ref{prop:weight}.

\begin{corollary}[Besov data and self-similar solutions]\label{cor:data}
Suppose that the hypotheses of Theorem \ref{thm:global} hold, and let
\begin{equation}\label{eq:sigma0}
\sigma_0:=\frac n{r'}-\frac{2-\gamma}{p-1}.
\end{equation}
\begin{enumerate}[(i)]
\item There exists $C>0$ such that, for all $(u_0,u_1)\in\B{\sigma_0}{r'}\infty\times\B{\sigma_0-1}{r'}\infty$,
\begin{equation}\label{eq:data-embedding}
\|(u_0,u_1)\|_{\Ic}\le C\bigl(\|u_0\|_{\B{\sigma_0}{r'}\infty}+\|u_1\|_{\B{\sigma_0-1}{r'}\infty}\bigr).
\end{equation}
Consequently, the conclusions of Theorem \ref{thm:global} (i) hold whenever $\|u_0\|_{\B{\sigma_0}{r'}\infty}+\|u_1\|_{\B{\sigma_0-1}{r'}\infty}$ is small enough.
\item Let $\Omega_0,\Omega_1\in C^\infty(\mathbb S^{n-1})$, and define $u_0,u_1$ by \eqref{eq:homogeneous}. Then $(u_0,u_1)\in\B{\sigma_0}{r'}\infty\times\B{\sigma_0-1}{r'}\infty$, and there exists $\varepsilon_0>0$ such that, for $0<|\epsilon|\le\varepsilon_0$, the solution of \eqref{eq:integral} with data $(\epsilon u_0,\epsilon u_1)$ given by Theorem \ref{thm:global} is self-similar.
\end{enumerate}
\end{corollary}

\begin{proof}
(i) By \eqref{eq:H1}, $\eta_r>\frac1p$, and \eqref{eq:zeta} gives $(p-1)\zeta=1-\eta_r<(p-1)\eta_r$, that is, $\zeta<\eta_r$. We apply Lemma \ref{lem:group} with $\sigma=\infty$ and $\theta:=\zeta/\eta_r\in(0,1)$, so that $\theta\eta_r=\zeta$. Since $n(1-\frac2r)-2\beta_r=(n-1)(\frac12-\frac1r)=\eta_r$, we have
\[
(1-\theta)n\Bigl(1-\frac2r\Bigr)+2\theta\beta_r=n\Bigl(1-\frac2r\Bigr)-\theta\eta_r=n\Bigl(1-\frac2r\Bigr)-\zeta,
\]
and \eqref{eq:zeta} gives
\[
s+n\Bigl(1-\frac2r\Bigr)-\zeta=n-\frac nr-\frac{2-\gamma}{p-1}=\sigma_0 .
\]
Therefore \eqref{eq:group-dK} and \eqref{eq:group-K} yield
\[
t^\zeta\|\dK(t)u_0\|_{\B sr\infty}\le C\|u_0\|_{\B{\sigma_0}{r'}\infty},\qquad
t^\zeta\|K(t)u_1\|_{\B sr\infty}\le C\|u_1\|_{\B{\sigma_0-1}{r'}\infty}
\]
for all $t>0$, which is \eqref{eq:data-embedding}. Observe that $\sigma_0-1<\sigma_0<\frac n{r'}$, so the data spaces consist of distributions in $\Sch'_h(\Rn)$. The last assertion of (i) follows from Theorem \ref{thm:global} (i).

(ii) Let $\alpha_0:=\frac{2-\gamma}{p-1}$. Since $p>2$, $0<\alpha_0<2\le n$. We also have $\alpha_0+1<n$. This is clear for $n\ge3$. For $n=2$, Remark \ref{rem:hypotheses} (b) gives $p>\frac3{1-\gamma}>3$, hence $\alpha_0<1$. Proposition \ref{prop:weight} with $m=r'$ therefore gives $u_0\in\B{\frac n{r'}-\alpha_0}{r'}\infty=\B{\sigma_0}{r'}\infty$ and $u_1\in\B{\frac n{r'}-\alpha_0-1}{r'}\infty=\B{\sigma_0-1}{r'}\infty$. By (i), $\|(\epsilon u_0,\epsilon u_1)\|_{\Ic}\le C|\epsilon|\bigl(\|u_0\|_{\B{\sigma_0}{r'}\infty}+\|u_1\|_{\B{\sigma_0-1}{r'}\infty}\bigr)\le\varepsilon$ for $|\epsilon|$ small. Theorem \ref{thm:global} (ii) then applies, because $\epsilon u_0$ and $\epsilon u_1$ are homogeneous of the right degrees.
\end{proof}

\begin{remark}\label{rem:final}
Compared with methods based on Strichartz estimates, the time-weighted framework only uses the dispersive estimate \eqref{eq:group}, and the solution is measured on the Besov space $\B{s_r}r\infty$. There is a price to pay. The time-integrability conditions $\eta_r<1$ and $p\zeta_r<1$ restrict the range of $r$. Moreover, the composition estimate has to produce the $L^{r'}$-based regularity $\sigma_r$ from the $L^r$-based regularity $s_r$, which forces us to estimate the coefficients $M_k(u)$ in $L^{r/(r-2)}$. This is where condition \eqref{eq:H2} comes from, and \eqref{eq:H2} in turn forces $\gamma<1$.
\end{remark}

\section{Scattering and asymptotic behavior}\label{sec:asymptotic}

We now describe the behavior of the solutions of Theorem \ref{thm:global} for large times. In this section the hypotheses of Theorem \ref{thm:global} hold, $s=s_r$, $\zeta=\zeta_r$, and $\varepsilon$ is the constant of Theorem \ref{thm:global}. For $\mu\ge0$ and $w:(0,\infty)\to\B sr\infty$ we use the notation
\[
\|w\|_{E_{s,\mu}}:=\sup_{t>0}t^\mu\|w(t)\|_{\B sr\infty},
\]
which agrees with Definition \ref{def:E} when $\mu=\zeta$. For a solution $u$ of \eqref{eq:integral} we write $$G(\tau):=a\,h_\gamma F(u(\tau)),$$
so that $B(u)(t)=\int_0^tK(t-\tau)G(\tau)d\tau$.

We first deduce two consequences of the previous sections. From Lemma \ref{lem:group-potential} with $\theta=1$, $\sigma=\infty$ and $\varsigma=\sigma_r$, and from \eqref{eq:unweighted-diff}, we obtain a constant $C_*>0$, depending on $n,\gamma,p,r$, such that
\begin{equation}\label{eq:Cstar}
\bigl\|K(t)\bigl(h_\gamma(F(w_1)-F(w_2))\bigr)\bigr\|_{\B sr\infty}\le C_*|t|^{-\eta_r}\bigl(\|w_1\|_{\B sr\infty}+\|w_2\|_{\B sr\infty}\bigr)^{p-1}\|w_1-w_2\|_{\B sr\infty}
\end{equation}
for all $t\ne0$ and $w_1,w_2\in\B sr\infty$. In particular, since $F(0)=0$,
\[
\|K(t)(h_\gamma F(w))\|_{\B sr\infty}\le C_*|t|^{-\eta_r}\|w\|_{\B sr\infty}^p,\qquad t\ne0,\ w\in\B sr\infty .
\]
Furthermore, \eqref{eq:beta-eta} and \eqref{eq:zeta} give
\begin{equation}\label{eq:exponent-relations}
\eta_r+(p-1)\zeta=1,\qquad \eta_r+p\zeta=1+\zeta,\qquad n\Bigl(1-\frac2r\Bigr)-2\beta_r=\eta_r .
\end{equation}

\subsection{Scattering}

We start with the scattering property. Every solution given by Theorem \ref{thm:global} (i) approaches a free wave $u^+$ as $t\to\infty$. The data of $u^+$ are the limits, as $t\to\infty$, of the data at time zero of the free waves that coincide with $u$ at time $t$. The next theorem makes this precise, with the rate $t^{-\zeta}$ for the convergence of these data. Its proof relies on the block estimates for the wave group from Lemma \ref{lem:group} and on \eqref{eq:Cstar}.

\begin{theorem}[Scattering]\label{thm:scattering}
Suppose that the hypotheses of Theorem \ref{thm:global} hold, and let $u$ be the solution given by Theorem \ref{thm:global} (i). Then there exists $(u_0^+,u_1^+)\in\Ic$ such that $\|(u_0^+,u_1^+)\|_{\Ic}\le\|u\|_{\E}+C|a|\,\|u\|_{\E}^p$ and the free wave $u^+(t):=\dK(t)u_0^++K(t)u_1^+$ satisfies
\begin{equation}\label{eq:scattering}
t^\zeta\|u(t)-u^+(t)\|_{\B sr\infty}\le C|a|\Bigl(\sup_{\tau\ge t}\tau^\zeta\|u(\tau)\|_{\B sr\infty}\Bigr)^p,\qquad t>0 .
\end{equation}
Moreover, for $t>0$ the pair
\[
V_0(t):=u_0-\int_0^tK(\tau)G(\tau)\,d\tau,\qquad V_1(t):=u_1+\int_0^t\dK(\tau)G(\tau)\,d\tau
\]
satisfies $u(t)=\dK(t)V_0(t)+K(t)V_1(t)$ and
\begin{equation}\label{eq:scattering-profile}
\|V_0(t)-u_0^+\|_{\B sr\infty}+\|V_1(t)-u_1^+\|_{\B{s-1}r\infty}\le C|a|\,\|u\|_{\E}^p\,t^{-\zeta}.
\end{equation}
In particular, $(V_0(t),V_1(t))\to(u_0^+,u_1^+)$ in $\B sr\infty\times\B{s-1}r\infty$ as $t\to\infty$. The scattering data are given by
\begin{equation}\label{eq:scattering-data}
u_0^+=u_0-\int_0^\infty K(\tau)G(\tau)\,d\tau,\qquad u_1^+=u_1+\int_0^\infty\dK(\tau)G(\tau)\,d\tau,
\end{equation}
and the constant $C$ depends only on $n,\gamma,p,r$.
\end{theorem}

\begin{proof}
Proposition \ref{prop:weighted-nonlinear} and the definition of $\E$ give
\begin{equation}\label{eq:G-bound}
\|G(\tau)\|_{\B{s-1+2\beta_r}{r'}\infty}\le C|a|\,\|u(\tau)\|_{\B sr\infty}^p\le C|a|\,\tau^{-p\zeta}\|u\|_{\E}^p,\qquad\tau>0 .
\end{equation}
Consequently, $\sum_{|k-j|\le2}\|\Dj kG(\tau)\|_{L^{r'}}\le C|a|2^{-j(s-1+2\beta_r)}\tau^{-p\zeta}\|u\|_{\E}^p$. Using \eqref{eq:group-block} with $\theta=1$ and with $\theta=0$, together with \eqref{eq:exponent-relations}, we obtain, for all $j\in\Z$ and $\tau>0$,
\begin{align}
\|\Dj jK(\tau)G(\tau)\|_{L^r}&\le C|a|\,\|u\|_{\E}^p\min\bigl\{\tau^{-1-\zeta}2^{-js},\ \tau^{-p\zeta}2^{-j(s-\eta_r)}\bigr\},\label{eq:block-K}\\
\|\Dj j\dK(\tau)G(\tau)\|_{L^r}&\le C|a|\,\|u\|_{\E}^p\min\bigl\{\tau^{-1-\zeta}2^{-j(s-1)},\ \tau^{-p\zeta}2^{-j(s-1-\eta_r)}\bigr\}.\label{eq:block-dK}
\end{align}
For $\theta=1$, the exponent of $2^j$ in \eqref{eq:group-block} is $-1+2\beta_r$, the power of $\tau$ is $\tau^{-\eta_r}$, and $\eta_r+p\zeta=1+\zeta$. For $\theta=0$, the exponent of $2^j$ is $-1+n(1-\frac2r)$, and $n(1-\frac2r)-2\beta_r=\eta_r$. Since $p\zeta<1$ and $\zeta>0$, \eqref{eq:block-K} implies
\[
\begin{aligned}
2^{j(s-\eta_r)}\int_0^1\|\Dj jK(\tau)G(\tau)\|_{L^r}d\tau&\le\frac{C|a|\,\|u\|_{\E}^p}{1-p\zeta},\\
2^{js}\int_t^\infty\|\Dj jK(\tau)G(\tau)\|_{L^r}d\tau&\le\frac{C|a|\,\|u\|_{\E}^p}{\zeta}\,t^{-\zeta},
\end{aligned}
\]
for all $j\in\Z$ and $t>0$. The same argument applied to \eqref{eq:block-dK} gives the analogous bounds for $\dK$, with $s-1-\eta_r$ and $s-1$ instead of $s-\eta_r$ and $s$. Then
\begin{equation}\label{eq:tails}
\Bigl\|\int_t^\infty K(\tau)G(\tau)d\tau\Bigr\|_{\B sr\infty}+\Bigl\|\int_t^\infty\dK(\tau)G(\tau)d\tau\Bigr\|_{\B{s-1}r\infty}\le C|a|\,\|u\|_{\E}^p\,t^{-\zeta},\qquad t>0 .
\end{equation}
As $s-1-\eta_r<s-\eta_r<s<\frac nr$, these Besov spaces are Banach spaces contained in $\Sch'_h(\Rn)$. Hence $(u_0^+,u_1^+)$ is a well-defined element of $\Sch'_h(\Rn)\times\Sch'_h(\Rn)$, and so are $V_0(t)$ and $V_1(t)$ for each $t>0$. Moreover,
\[
V_0(t)-u_0^+=\int_t^\infty K(\tau)G(\tau)\,d\tau,\qquad V_1(t)-u_1^+=-\int_t^\infty\dK(\tau)G(\tau)\,d\tau .
\]
Thus \eqref{eq:tails} is exactly \eqref{eq:scattering-profile}.

The addition formula $\sin((t-\tau)|\xi|)=\sin(t|\xi|)\cos(\tau|\xi|)-\cos(t|\xi|)\sin(\tau|\xi|)$ gives
\begin{equation}\label{eq:K-identity}
K(t-\tau)=K(t)\dK(\tau)-\dK(t)K(\tau),\qquad t,\tau\in\R .
\end{equation}
The operators $\Dj jK(t)$ and $\Dj j\dK(t)$ are bounded on $L^r(\Rn)$, so they commute with the integrals above. Together with \eqref{eq:K-identity}, this gives, for $j\in\Z$ and $t>0$,
\[
\Dj j\bigl(\dK(t)V_0(t)+K(t)V_1(t)\bigr)=\Dj ju_L(t)+\int_0^t\Dj jK(t-\tau)G(\tau)\,d\tau=\Dj ju(t)
\]
and
\[
\Dj ju^+(t)=\Dj ju_L(t)+\int_0^\infty\Dj jK(t-\tau)G(\tau)\,d\tau .
\]
The last integral converges absolutely in $L^r$ by \eqref{eq:group-block} and \eqref{eq:G-bound}. In particular, $u(t)=\dK(t)V_0(t)+K(t)V_1(t)$, and
\begin{equation}\label{eq:u-minus-uplus}
\Dj j\bigl(u(t)-u^+(t)\bigr)=-\int_t^\infty\Dj jK(t-\tau)G(\tau)\,d\tau,\qquad j\in\Z,\ t>0 .
\end{equation}
Let $N(t):=\sup_{\tau\ge t}\tau^\zeta\|u(\tau)\|_{\B sr\infty}$. Using \eqref{eq:u-minus-uplus} and \eqref{eq:Cstar}, we find
\[
\begin{aligned}
\|u(t)-u^+(t)\|_{\B sr\infty}&\le\int_t^\infty\|K(t-\tau)G(\tau)\|_{\B sr\infty}d\tau\le C_*|a|\int_t^\infty(\tau-t)^{-\eta_r}\|u(\tau)\|_{\B sr\infty}^pd\tau\\
&\le C_*|a|\,N(t)^p\int_t^\infty(\tau-t)^{-\eta_r}\tau^{-p\zeta}d\tau
=C_*|a|\,N(t)^p\,t^{-\zeta}\int_1^\infty(z-1)^{-\eta_r}z^{-p\zeta}dz .
\end{aligned}
\]
Here $1-\eta_r-p\zeta=-\zeta$ by \eqref{eq:exponent-relations}, and the last integral is finite because $\eta_r<1$ and $\eta_r+p\zeta=1+\zeta>1$. This proves \eqref{eq:scattering}. Finally, $N(t)\le\|u\|_{\E}$, and \eqref{eq:scattering} gives $t^\zeta\|u^+(t)\|_{\B sr\infty}\le\|u\|_{\E}+C|a|\,\|u\|_{\E}^p$ for every $t>0$. Hence $(u_0^+,u_1^+)\in\Ic$ with the stated bound.
\end{proof}

\subsection{Decay}

The solutions of Theorem \ref{thm:global} decay like $t^{-\zeta}$ in $\B sr\infty$, and this rate is dictated by the scaling. The next theorem shows that a solution keeps the decay rate $t^{-\mu}$ of its free evolution $u_L$ for every $\mu\in[0,\eta_r)$, provided the data are small enough. Since $\zeta<\eta_r$, this includes rates faster than $t^{-\zeta}$. When $\mu>\zeta$, the difference $u(t)-u^+(t)$ decays even faster than $t^{-\mu}$. The proof repeats the fixed point argument of Theorem \ref{thm:global} in a set where the norm of $E_{s,\mu}$ is also controlled. To this end, 
for $\mu\in[0,\eta_r)$ we define
\[
\mathbf B_\mu:=\int_0^1(1-z)^{-\eta_r}z^{-(1-\eta_r)-\mu}\,dz,
\]
which is finite because $\eta_r<1$ and $(1-\eta_r)+\mu<1$. We choose $\varepsilon_\mu\in(0,\varepsilon]$ such that
\begin{equation}\label{eq:eps-mu}
C_*|a|\,(4\varepsilon_\mu)^{p-1}\mathbf B_\mu\le\frac12 .
\end{equation}
Let $w_1,w_2\in\E$ with $\|w_1\|_{\E},\|w_2\|_{\E}\le2\varepsilon_\mu$, so that $(\|w_1(\tau)\|_{\B sr\infty}+\|w_2(\tau)\|_{\B sr\infty})^{p-1}\le(4\varepsilon_\mu)^{p-1}\tau^{-(p-1)\zeta}$ for $\tau>0$. By \eqref{eq:Cstar} and the substitution $\tau=tz$ we have, for $t>0$,
\begin{equation}\label{eq:B-mu}
\begin{aligned}
&t^\mu\|B(w_1)(t)-B(w_2)(t)\|_{\B sr\infty}\\
&\quad\le C_*|a|\,(4\varepsilon_\mu)^{p-1}t^\mu\int_0^t(t-\tau)^{-\eta_r}\tau^{-(p-1)\zeta}\|w_1(\tau)-w_2(\tau)\|_{\B sr\infty}d\tau\\
&\quad=C_*|a|\,(4\varepsilon_\mu)^{p-1}\int_0^1(1-z)^{-\eta_r}z^{-(1-\eta_r)-\mu}\,(tz)^\mu\|w_1(tz)-w_2(tz)\|_{\B sr\infty}dz .
\end{aligned}
\end{equation}
Here, we used $(p-1)\zeta=1-\eta_r$, which shows that the power of $t$ is $\mu+1-\eta_r-(1-\eta_r)-\mu=0$. In particular, \eqref{eq:eps-mu} gives
\begin{equation}\label{eq:B-mu-2}
\|B(w_1)-B(w_2)\|_{E_{s,\mu}}\le\frac12\|w_1-w_2\|_{E_{s,\mu}} .
\end{equation}

\begin{theorem}[Decay]\label{thm:decay}
Suppose that the hypotheses of Theorem \ref{thm:global} hold, and let $\mu\in[0,\eta_r)$. Then there exists $\varepsilon_\mu\in(0,\varepsilon]$ such that, if $(u_0,u_1)\in\Ic$ satisfies $\|(u_0,u_1)\|_{\Ic}\le\varepsilon_\mu$ and $\|u_L\|_{E_{s,\mu}}<\infty$, the solution $u$ given by Theorem \ref{thm:global} has the following properties.
\begin{enumerate}[(i)]
\item $\|u\|_{E_{s,\mu}}\le2\|u_L\|_{E_{s,\mu}}$.
\item If $\mu>\zeta$, then the free wave $u^+$ of Theorem \ref{thm:scattering} satisfies
\begin{equation}\label{eq:decay-scattering}
\|u(t)-u^+(t)\|_{\B sr\infty}\le C|a|\,\|u\|_{E_{s,\mu}}^p\,t^{-\mu-(p-1)(\mu-\zeta)},\qquad t>0 .
\end{equation}
\end{enumerate}
\end{theorem}

\begin{proof}
(i) Let $(u_0,u_1)\in\Ic$ with $\|(u_0,u_1)\|_{\Ic}\le\varepsilon_\mu$ and $\|u_L\|_{E_{s,\mu}}<\infty$, and set
\[
Y:=\bigl\{w\in\E:\ \|w\|_{\E}\le2\varepsilon_\mu,\ \|w\|_{E_{s,\mu}}\le2\|u_L\|_{E_{s,\mu}}\bigr\}.
\]
The set $Y$ contains $u_L$. It is closed in $\E$, because convergence in $\E$ implies convergence of $w(t)$ in $\B sr\infty$ for each $t>0$. Let $\Phi(w)=u_L+B(w)$ be the map used in the proof of Theorem \ref{thm:global}. Since $\varepsilon_\mu\le\varepsilon$, that proof shows that $\|\Phi(w)\|_{\E}\le2\varepsilon_\mu$ whenever $\|w\|_{\E}\le2\varepsilon_\mu$, and that $\Phi$ is a contraction on $\overline B(0,2\varepsilon)$. For $w\in Y$, \eqref{eq:B-mu-2} with $w_1=w$, $w_2=0$ and $B(0)=0$ gives
\[
\|\Phi(w)\|_{E_{s,\mu}}\le\|u_L\|_{E_{s,\mu}}+\frac12\|w\|_{E_{s,\mu}}\le2\|u_L\|_{E_{s,\mu}} .
\]
Hence $\Phi$ maps $Y$ into itself, and it is a contraction on $Y$ for the norm of $\E$. It therefore has a fixed point in $Y$. By uniqueness in $\overline B(0,2\varepsilon)$, this fixed point is the solution $u$, and so $\|u\|_{E_{s,\mu}}\le2\|u_L\|_{E_{s,\mu}}$.

(ii) Let $\mu>\zeta$. For $t>0$, \eqref{eq:u-minus-uplus} and \eqref{eq:Cstar} give
\[
\begin{aligned}
\|u(t)-u^+(t)\|_{\B sr\infty}&\le C_*|a|\,\|u\|_{E_{s,\mu}}^p\int_t^\infty(\tau-t)^{-\eta_r}\tau^{-p\mu}d\tau\\
&=C_*|a|\,\|u\|_{E_{s,\mu}}^p\,t^{1-\eta_r-p\mu}\int_1^\infty(z-1)^{-\eta_r}z^{-p\mu}dz .
\end{aligned}
\]
The last integral is finite, since $\eta_r<1$ and $\eta_r+p\mu>\eta_r+p\zeta=1+\zeta>1$. As $1-\eta_r=(p-1)\zeta$, we have $1-\eta_r-p\mu=-\mu-(p-1)(\mu-\zeta)$, which proves \eqref{eq:decay-scattering}.
\end{proof}

\begin{remark}\label{rem:decay}
By \eqref{eq:H1} and \eqref{eq:zeta}, $\zeta<\eta_r$, so Theorem \ref{thm:decay} includes decay rates faster than $t^{-\zeta}$. For $\mu\in[0,\eta_r]$, the condition $\|u_L\|_{E_{s,\mu}}<\infty$ holds whenever
\[
(u_0,u_1)\in\B{\sigma_0+\zeta-\mu}{r'}\infty\times\B{\sigma_0-1+\zeta-\mu}{r'}\infty .
\]
Indeed, Lemma \ref{lem:group} with $\theta=\mu/\eta_r$ gives $t^\mu\|u_L(t)\|_{\B sr\infty}\le C(\|u_0\|_{\B{\sigma_0+\zeta-\mu}{r'}\infty}+\|u_1\|_{\B{\sigma_0-1+\zeta-\mu}{r'}\infty})$, because $s+(1-\theta)n(1-\frac2r)+2\theta\beta_r=s+n(1-\frac2r)-\mu=\sigma_0+\zeta-\mu$ (see the proof of Corollary \ref{cor:data}). For $\mu=\zeta$ this is \eqref{eq:data-embedding}.
\end{remark}

\subsection{Asymptotic stability}

In what follows, we `compare two solutions. The next theorem shows that the large time behavior of a solution is determined by that of its free evolution. More precisely, if the difference of two free evolutions decays faster than $t^{-\mu}$, then so does the difference of the corresponding solutions. The proof is based on the integral inequality \eqref{eq:B-mu} and on the Fatou lemma.

\begin{theorem}[Asymptotic stability]\label{thm:stability}
Suppose that the hypotheses of Theorem \ref{thm:global} hold, let $\mu\in[0,\eta_r)$, and let $\varepsilon_\mu$ be given by Theorem \ref{thm:decay}. Let $(u_0,u_1),(v_0,v_1)\in\Ic$ be such that $\|(u_0,u_1)\|_{\Ic}\le\varepsilon_\mu$, $\|(v_0,v_1)\|_{\Ic}\le\varepsilon_\mu$ and $\|u_L\|_{E_{s,\mu}}+\|v_L\|_{E_{s,\mu}}<\infty$, where $v_L(t):=\dK(t)v_0+K(t)v_1$, and let $u$ and $v$ be the corresponding solutions given by Theorem \ref{thm:global}. If
\begin{equation}\label{eq:stability-hyp}
\lim_{t\to\infty}t^\mu\|u_L(t)-v_L(t)\|_{\B sr\infty}=0,
\end{equation}
then
\begin{equation}\label{eq:stability}
\lim_{t\to\infty}t^\mu\|u(t)-v(t)\|_{\B sr\infty}=0 .
\end{equation}
For $\mu=\zeta$, the condition $\|u_L\|_{E_{s,\mu}}+\|v_L\|_{E_{s,\mu}}<\infty$ holds automatically.
\end{theorem}

\begin{proof}
By Theorem \ref{thm:decay} (i), $D:=\|u-v\|_{E_{s,\mu}}\le\|u\|_{E_{s,\mu}}+\|v\|_{E_{s,\mu}}<\infty$. In addition, $\|u\|_{\E},\|v\|_{\E}\le2\varepsilon_\mu$, since $u$ and $v$ lie in the corresponding sets $Y$ from the proof of Theorem \ref{thm:decay} (i). Writing $u-v=u_L-v_L+B(u)-B(v)$ and using \eqref{eq:B-mu}, we get, for $t>0$,
\begin{equation}\label{eq:stability-ineq}
t^\mu\|u(t)-v(t)\|_{\B sr\infty}\le t^\mu\|u_L(t)-v_L(t)\|_{\B sr\infty}+C_*|a|\,(4\varepsilon_\mu)^{p-1}\int_0^1(1-z)^{-\eta_r}z^{-(1-\eta_r)-\mu}\,g(tz)\,dz,
\end{equation}
where $g(\tau):=\tau^\mu\|u(\tau)-v(\tau)\|_{\B sr\infty}\le D$. Let $A:=\limsup_{t\to\infty}g(t)\in[0,D]$. For each $z\in(0,1)$ we have $\limsup_{t\to\infty}g(tz)=A$. Moreover, the integrand in \eqref{eq:stability-ineq} is dominated by the integrable function $D(1-z)^{-\eta_r}z^{-(1-\eta_r)-\mu}$. We take the $\limsup$ as $t\to\infty$ in \eqref{eq:stability-ineq}, use \eqref{eq:stability-hyp}, and apply the Fatou lemma to the nonnegative functions $D(1-z)^{-\eta_r}z^{-(1-\eta_r)-\mu}-(1-z)^{-\eta_r}z^{-(1-\eta_r)-\mu}g(tz)$. By applying \eqref{eq:eps-mu}, we can conclude that
\[
A\le0+C_*|a|\,(4\varepsilon_\mu)^{p-1}\mathbf B_\mu A\le\frac12A .
\]
Hence $A=0$, which is \eqref{eq:stability}. If $\mu=\zeta$, then $\|u_L\|_{E_{s,\zeta}}=\|(u_0,u_1)\|_{\Ic}$ and $\|v_L\|_{E_{s,\zeta}}=\|(v_0,v_1)\|_{\Ic}$, which are finite.
\end{proof}

\subsection{Asymptotic self-similarity}

As a consequence, small perturbations of homogeneous data produce solutions that approach self-similar solutions for large times. The perturbations are taken in Besov spaces with finite third index $q$, which guarantees that their free evolution decays faster than $t^{-\zeta}$. The following result then follows from Theorem \ref{thm:stability} with $\mu=\zeta$, applied to the solutions with and without the perturbation. To this end, we set $\alpha_0:=\frac{2-\gamma}{p-1}$.

\begin{corollary}[Asymptotic self-similarity]\label{cor:asymptotic-ss}
Suppose that the hypotheses of Theorem \ref{thm:global} hold. Let $\Omega_0,\Omega_1\in C^\infty(\mathbb S^{n-1})$, let $\varphi_0,\varphi_1$ be the homogeneous functions on the right-hand sides of \eqref{eq:homogeneous}, and let $1\le q<\infty$. Then there exists $\delta>0$ with the following property. If $\epsilon\in\R$ and $(w_0,w_1)\in\B{\sigma_0}{r'}q\times\B{\sigma_0-1}{r'}q$ satisfy
\[
|\epsilon|\bigl(\|\varphi_0\|_{\B{\sigma_0}{r'}\infty}+\|\varphi_1\|_{\B{\sigma_0-1}{r'}\infty}\bigr)+\|w_0\|_{\B{\sigma_0}{r'}q}+\|w_1\|_{\B{\sigma_0-1}{r'}q}\le\delta ,
\]
then the solution $v$ of \eqref{eq:integral} with data $(\epsilon\varphi_0,\epsilon\varphi_1)$ is self-similar, and the solution $u$ with data $(\epsilon\varphi_0+w_0,\epsilon\varphi_1+w_1)$ satisfies
\begin{equation}\label{eq:asymptotic-ss}
\lim_{t\to\infty}t^\zeta\|u(t)-v(t)\|_{\B sr\infty}=0 .
\end{equation}
Moreover, $v(t,x)=t^{-\alpha_0}V(x/t)$ with $V:=v(1)$, and, for every $t_0>0$,
\begin{equation}\label{eq:asymptotic-profile}
\lim_{\lambda\to\infty}\ \sup_{t\ge t_0}\,t^\zeta\bigl\|\lambda^{\alpha_0}u(\lambda t,\lambda\,\cdot)-v(t)\bigr\|_{\B sr\infty}=0 .
\end{equation}
\end{corollary}

\begin{proof}
Let $C_{\mathcal I}$ be the constant in \eqref{eq:data-embedding}, let $\varepsilon_\zeta$ be given by Theorem \ref{thm:decay} with $\mu=\zeta$, and set $\delta:=\varepsilon_\zeta/C_{\mathcal I}$. By Corollary \ref{cor:data} (ii), $\varphi_0\in\B{\sigma_0}{r'}\infty$ and $\varphi_1\in\B{\sigma_0-1}{r'}\infty$. Since $\|\cdot\|_{\B{\sigma}{r'}\infty}\le\|\cdot\|_{\B{\sigma}{r'}q}$, \eqref{eq:data-embedding} gives
\[
\|(\epsilon\varphi_0,\epsilon\varphi_1)\|_{\Ic}\le\varepsilon_\zeta,\qquad\|(\epsilon\varphi_0+w_0,\epsilon\varphi_1+w_1)\|_{\Ic}\le\varepsilon_\zeta .
\]
As $\varepsilon_\zeta\le\varepsilon$ and $\epsilon\varphi_0$, $\epsilon\varphi_1$ are homogeneous of degrees $-\alpha_0$ and $-\alpha_0-1$, Theorem \ref{thm:global} (ii) shows that $v$ is self-similar. Moreover, $u_L-v_L=w_L$ with $w_L(t):=\dK(t)w_0+K(t)w_1$. By Theorem \ref{thm:stability} with $\mu=\zeta$, it is therefore enough to prove that $t^\zeta\|w_L(t)\|_{\B sr\infty}\to0$ as $t\to\infty$.

For $N\in\mathbb N$ let $w_i^N:=\sum_{|j|\le N}\Dj jw_i$, $i=0,1$, and $w_L^N(t):=\dK(t)w_0^N+K(t)w_1^N$. Note that $\Dj k(w_i-w_i^N)=0$ for $|k|\le N-2$, and that otherwise $\Dj k(w_i-w_i^N)$ is a sum of at most three blocks of $w_i$. Since $q<\infty$, this gives
\[
\begin{aligned}
\varrho_N&:=\|w_0-w_0^N\|_{\B{\sigma_0}{r'}\infty}+\|w_1-w_1^N\|_{\B{\sigma_0-1}{r'}\infty}\\
&\le C\Bigl(\sum_{|k|\ge N-3}\bigl(2^{k\sigma_0}\|\Dj kw_0\|_{L^{r'}}+2^{k(\sigma_0-1)}\|\Dj kw_1\|_{L^{r'}}\bigr)^q\Bigr)^{\frac1q}\to0
\end{aligned}
\]
as $N\to\infty$. By \eqref{eq:data-embedding}, $t^\zeta\|w_L(t)-w_L^N(t)\|_{\B sr\infty}\le C_{\mathcal I}\varrho_N$ for every $t>0$. On the other hand, $w_0^N$ and $w_1^N$ have only finitely many nonzero blocks. Hence $w_0^N\in\B{s+2\beta_r}{r'}\infty$ and $w_1^N\in\B{s-1+2\beta_r}{r'}\infty$, and Lemma \ref{lem:group} with $\theta=1$ gives
\[
t^\zeta\|w_L^N(t)\|_{\B sr\infty}\le Ct^{\zeta-\eta_r}\bigl(\|w_0^N\|_{\B{s+2\beta_r}{r'}\infty}+\|w_1^N\|_{\B{s-1+2\beta_r}{r'}\infty}\bigr)\to0\qquad\text{as }t\to\infty,
\]
because $\zeta<\eta_r$. Consequently, $\limsup_{t\to\infty}t^\zeta\|w_L(t)\|_{\B sr\infty}\le C_{\mathcal I}\varrho_N$ for every $N$. Letting $N\to\infty$, we obtain \eqref{eq:asymptotic-ss}.

Finally, $v=v_\lambda$ for every $\lambda>0$, and the choice $\lambda=1/t$ gives $v(t,x)=t^{-\alpha_0}v(1,x/t)$. Moreover, $\lambda^{\alpha_0}u(\lambda t,\lambda\cdot)-v(t)=(u-v)_\lambda(t)$, and, as in the proof of Theorem \ref{thm:global} (ii), $\|(u-v)_\lambda(t)\|_{\B sr\infty}\le C\lambda^{\alpha_0+s-\frac nr}\|(u-v)(\lambda t)\|_{\B sr\infty}=C\lambda^\zeta\|(u-v)(\lambda t)\|_{\B sr\infty}$. Therefore
\[
\sup_{t\ge t_0}t^\zeta\bigl\|\lambda^{\alpha_0}u(\lambda t,\lambda\,\cdot)-v(t)\bigr\|_{\B sr\infty}\le C\sup_{\tau\ge\lambda t_0}\tau^\zeta\|u(\tau)-v(\tau)\|_{\B sr\infty}\to0\qquad\text{as }\lambda\to\infty
\]
by \eqref{eq:asymptotic-ss}, which is \eqref{eq:asymptotic-profile}.
\end{proof}

\begin{remark}\label{rem:asymptotic-ss}
Since \eqref{eq:HHW} is invariant under \eqref{eq:scaling}, the attracting self-similar solution in Corollary \ref{cor:asymptotic-ss} solves the same equation, as in Cazenave and Weissler \cite{CW98}.  The assumption $q<\infty$ in Corollary \ref{cor:asymptotic-ss} may be replaced by $(w_0,w_1)\in(\B{\sigma_0}{r'}\infty\cap\B{\sigma_0+\zeta-\mu}{r'}\infty)\times(\B{\sigma_0-1}{r'}\infty\cap\B{\sigma_0-1+\zeta-\mu}{r'}\infty)$ for some $\mu\in(\zeta,\eta_r]$. Indeed, Remark \ref{rem:decay} then gives $t^\zeta\|\dK(t)w_0+K(t)w_1\|_{\B sr\infty}\le Ct^{\zeta-\mu}\to0$.
\end{remark}

\appendix
\section*{ Appendix}

\section{Paradifferential calculus in Besov spaces}\label{subsec:paradiff}

We use the notation of Bony's paradifferential calculus \cite{Bon81}, see \eqref{eq:LP}. For $f,g\in\Sch'_h(\Rn)$ the pointwise product $fg$ has the formal decomposition
\begin{equation}\label{eq:Bony}
fg=\sum_{j\in\Z}\Sj{j-3}f\,\Dj jg+\sum_{j\in\Z}\Sj{j-3}g\,\Dj jf+\sum_{j\in\Z}\tDj jf\,\Dj jg=:\dot T_fg+\dot T_gf+\dot R(f,g).
\end{equation}
With our choice of $\phi$, one checks easily that
\begin{equation}\label{eq:Bony-supp}
\Dj l(\dot T_fg)=\sum_{|j-l|\le4}\Dj l\bigl(\Sj{j-3}f\,\Dj jg\bigr),\qquad
\Dj l\bigl(\dot R(f,g)\bigr)=\sum_{j\ge l-5}\Dj l\bigl(\tDj jf\,\Dj jg\bigr).
\end{equation}
We refer to \cite[Chapter 2]{BCD11} for details.

We recall the product law of Abidi in homogeneous Besov spaces, which is used in Lemma \ref{lem:potential}. In this law one factor is only required to belong to a Besov space with third index $\infty$.

\begin{proposition}[{\cite[Proposition 2.2]{Abi07}}]\label{prop:Abidi}
Let $1\le m,\sigma\le\infty$ and $s_1,s_2\in\R$ satisfy
\begin{equation}\label{eq:Abidi-cond}
s_1<\frac nm,\qquad s_2<\frac nm,\qquad s_1+s_2+n\min\Bigl\{0,1-\frac2m\Bigr\}>0 .
\end{equation}
Then there is a constant $C>0$, depending on $n,m,s_1,s_2$ and $\sigma$, such that
\begin{equation}\label{eq:Abidi}
\|fg\|_{\B{s_1+s_2-\frac nm}m\sigma}\le C\|f\|_{\B{s_1}m\sigma}\|g\|_{\B{s_2}m\infty}.
\end{equation}
\end{proposition}

\end{document}